\documentclass[preprint,12pt]{elsarticle}

\usepackage{graphicx}
\usepackage{amssymb,amsfonts,amssymb}
\usepackage[margin=2.5cm]{geometry}

\usepackage{lineno}

\usepackage{graphics,graphicx}
\usepackage{epstopdf}

\usepackage{mathtools}
\usepackage{stmaryrd,dsfont}
\usepackage{bbm}
\usepackage{mathrsfs}
\usepackage{floatrow}
\usepackage{blindtext}
\usepackage{setspace}
\usepackage{pgfplots}
\usepackage{mathtools}
\usepackage{stmaryrd,dsfont}
\usepackage{bbm}
\usepackage{url}
\usepackage{wrapfig,bm}
\usepackage{amsmath,bm}
\usepackage[export]{adjustbox}
\usepackage{pgfplots}
\usepackage{tikz-3dplot}
\usepackage{booktabs}
\usepackage{color}
\usepackage{caption}
\usepackage[caption=false]{subfig}
\usepackage{float}
\floatstyle{plaintop}
\restylefloat{table}
\usepackage{tcolorbox}
\usepackage{tikz}
\usetikzlibrary{shapes.geometric,arrows}
\usepackage{fancyhdr,graphicx,amsmath,amssymb}
\usepackage[ruled,vlined]{algorithm2e}
\include{pythonlisting}
\biboptions{sort&compress}

\newcommand{\eref}[1]{Equation~(\ref{#1})}
\newcommand{\erefs}[1]{Equations~(\ref{#1})}
\newcommand{\fref}[1]{Figure~\ref{#1}}
\newcommand{\frefs}[1]{Figures~\ref{#1}}

\newcommand{\vm}[1]{\bm{\mathrm{#1}}} 
\newcommand{\bsym}[1]{\bm{#1}}
\newcommand{\cn}{\vm{n}}
\newcommand{\uu}{\vm{u}}

\newcommand{\bveps}{\bsym{\varepsilon}}
\newcommand{\bvsig}{\bsym{\sigma}}

\journal{Journal Name}

\begin{document}

\begin{frontmatter}


\title{Adaptive phase field method for quasi-static brittle fracture based on recovery based error indicator and quadtree decomposition}



\author{Hirshikesh, C Jansari, K Kannan, RK Annabattula}
\author{S Natarajan\corref{cor1}}
\address{Integrated Modeling and Simulation Lab, Department of Mechanical Engineering, Indian Institute of Technology Madras, Chennai-600036, India.}

\cortext[cor1]{Corresponding author}

\begin{abstract}
An adaptive phase field method is proposed for crack propagation in brittle materials under quasi-static loading. The adaptive refinement is based on the recovery type error indicator, which is combined with the quadtree decomposition. Such a decomposition leads to elements with hanging nodes. Thanks to the polygonal finite element method, the elements with hanging nodes are treated as polygonal elements and do not require any special treatment. The mean value coordinates are used to approximate the unknown field variables and a staggered solution scheme is adopted to compute the displacement and the phase field variable. A few standard benchmark problems are solved to show the efficiency of the proposed framework. It is seen that the proposed framework yields comparable results at a fraction of the computational cost when compared to standard approaches reported in the literature.
\end{abstract}

\begin{keyword}
crack propagation\sep phase field method \sep quadtree decomposition \sep quasi-static brittle fracture \sep recovery based error indicator \sep mean value coordinates


\end{keyword}

\end{frontmatter}


\section{Introduction}
\label{Section:introduction}

\begin{figure}[htpb]
\centering
\subfloat[]{\includegraphics[scale=0.6]{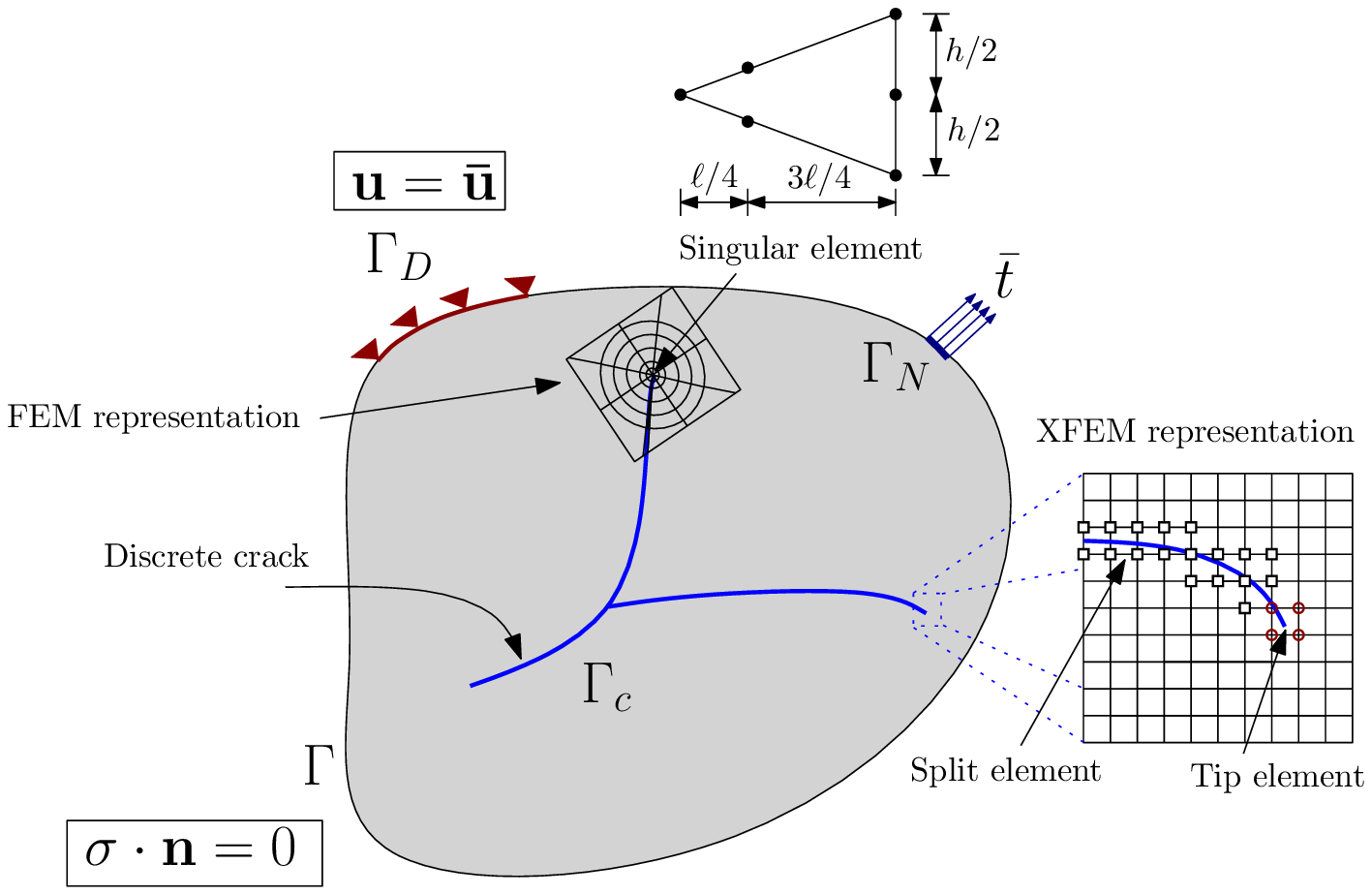}\label{discrete}} \hspace{5mm}
\subfloat[]{\includegraphics[scale = 0.55]{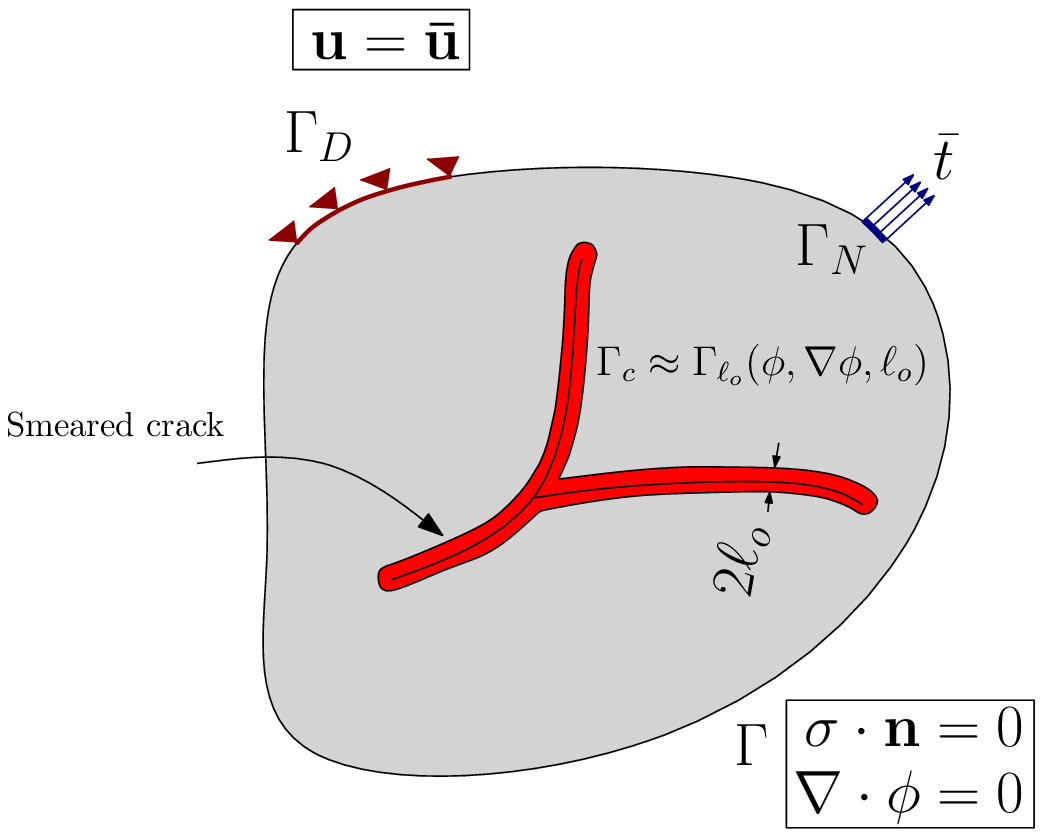} \label{diffuse} }
\caption{Schematic representation of a domain with a geometric discontinuity (a) discrete and (b) smeared/diffuse representation, where $\ell_o$ is the characteristic length scale.}
\label{fig:crackmodelapproaches}
\end{figure}

Two schools of thoughts that are commonly adopted to numerically simulate evolving discontinuities are: discrete approach and smeared or diffused approach. See \fref{fig:crackmodelapproaches} for a schematic illustration of different approaches. In case of discrete approach, the crack can either be represented explicitly or implicitly. For example, special elements such as singular elements and partition of unity methods, viz., PUFEM/GFEM/XFEM~\cite{Moes1999131, SANCHEZRIVADENEIRA2019876, Simone20061122, Gupta2015355,melenkbabuska1996,zhangbabuska2016} fall under this category. The discrete model with explicit representation, traces the crack morphology along a single or pre-defined finite number of interfaces. Frequent or continuous remeshing may be required to capture the changing morphology. On the other hand, implicit representation reduces the mesh burden by allowing the discontinuity to be independent of the background discretization. However, the success of the method relies on a priori knowledge of functions to capture the local behavior and a criteria to evolve the discontinuities \cite{BOUCHARD2003}. 
The smeared crack model does not require re-meshing and/or knowledge of additional function, however, necessitates enhancement of the element behavior to capture the kinematics, i.e., by modifying the stress strain relations. The introduction of variational approach to fracture~\cite{Francfort1998} has revolutionized the smeared crack approach. The method has been proven capable of modelling complex crack morphologies such as crack coalescence~\cite{MIEHE2016619}, crack branching~\cite{PAGGI2017145}, curvilinear crack paths~\cite{JEONG2018483}, to name a few. For a more comprehensive review of the phase field and its implementation, interested readers are referred to the work of \citet{Ambati2015,Wu2018,Bourdin2008,MANDAL201948}. 

Apart from this, the salient feature of the approach is that, existing finite element solvers can directly be used without any modification. However, the flexibility and robustness provided by the phase field method (PFM) comes with associated difficulties. The simplicity in handling complex crack morphologies comes: (a) at the cost of solving additional partial differential equation, viz., the phase field equation and (b) requires extremely fine meshes to accurately capture the crack topology~\cite{MSEKH2018287,BORDEN201277,Schlüter2014}. The latter, can be handled by local-refinement techniques, however, this requires the crack path to be known `{\it a priori}', which is often not the case. This aspect of the PFM has received considerable attention and has led to novel approaches for adaptive phase field model for crack propagation. They are based on the goal oriented error estimation adaptive grid strategy with the finite element method (FEM)~\cite{HEISTER2015466} and predictor-corrector strategies~\cite{MSEKH2018287, Burke2010988, AREIAS2016322, fabian2010689}. The predictor-corrector strategy coupled with $h$-adaptivity was employed to simulate thermo-mechanical cracks using the PFM~\cite{BADNAVA201831}. Zhou and Zhuang~\cite{ZHOU2018190} employed the PFM to simulate quasi-static crack propagation in rocks. Patil \textit{et al.,}~\cite{PATIL2018254} coupled the multiscale FEM with the hybrid PFM for brittle fracture problems. In this approach, the difference in the degrees of freedom between the coarse and the fine mesh is related by multiscale basis functions.

In this paper, by combining the adaptive refinement scheme based on recovery based error indicator with quadtree decomposition, we propose a novel adaptive phase field method for brittle fracture subject to quasi-static loading conditions. The salient features of the proposed work are: 
\begin{itemize}
    \item adaptive refinement scheme based on recovery based error indicator,
    \item quadtree decomposition technique which reduces the computational burden,
    \item hanging nodes due to quadtree decomposition are treated within the framework of the polygonal finite element method.
\end{itemize}

The rest of the paper is organized as follows: Section \ref{Section:numericalformulation} presents the governing equations for the phase field method and the corresponding weak form. The recovery based error indicator and the quadtree decomposition is discussed in Section \ref{Sec:Quadtree Meshing}. The robustness of the adaptive refinement technique is demonstrated with a few problems in Section \ref{Section:numericalexamples} and major conclusions are presented in the last section.

\section{Governing equations and the weak form}
\label{Section:numericalformulation}
\begin{figure}[htpb]
\centering
\includegraphics[scale = 0.7]{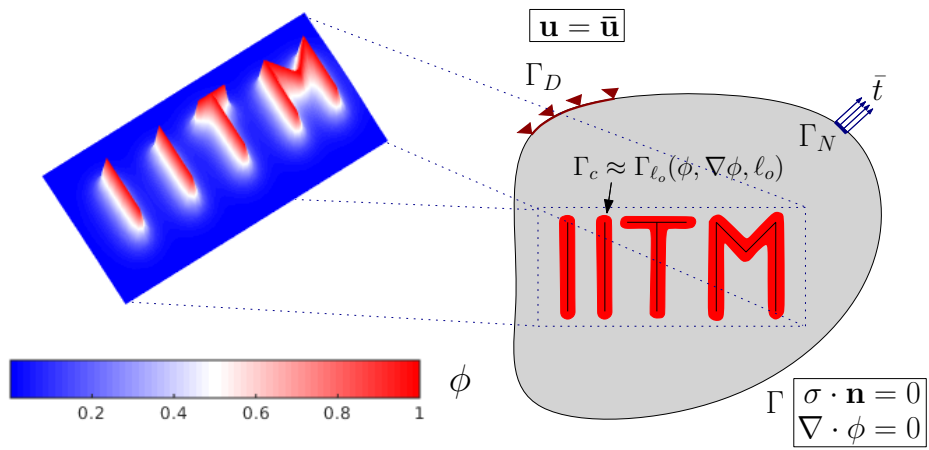} \label{diffuse2d}
\caption{Schematic representation of a domain with a geometric discontinuity in PFM framework, where $\phi \in [0, 1]$, is the phase-field variable, $\ell_o$ is the characteristic length scale.}
\label{domain_representation}
\end{figure}

Consider a linear elasto-static body with a discontinuity occupying the domain $\Omega \subset \mathbb{R}^d$, where $d=$ 2, 3 as shown in \fref{domain_representation}. The boundary ($\Gamma$) is considered to admit the decomposition with the outward normal $\cn$ into three disjoint sets, i.e., $\Gamma = \Gamma_{\rm D} \cup \Gamma_{\rm N} \cup \Gamma_{\rm c}$ and $\Gamma_{\rm D} \cap \Gamma_{\rm N} = \emptyset$, where on $\Gamma_{\rm D}$ and $\Gamma_{\rm N}$, Dirichlet boundary and Neumann boundary conditions are specified. The closure of the domain is $\overline{\Omega} \equiv \Omega \cup \Gamma$. In a discrete fracture mechanics context, the strong discontinuity, i.e., the crack is represented by a discontinuous surface ($\Gamma_{\rm c}$) as shown in \fref{discrete}. However, in the PFM, a diffuse field variable $\phi \in \left[0,1\right]$ is introduced to model the crack (see \fref{domain_representation}), with ${\rm max}(\phi)$ and ${\rm min}(\phi)$ representing the fully damaged and undamaged state, respectively.

In the absence of inertia and body forces, the coupled governing equations for linear elasto-statics undergoing small deformation is: find $(\uu,\phi): \Omega \rightarrow \mathbb{R}^d$ such that
\begin{subequations}
\begin{align}
\left[ (1-\phi)^2 + k_p \right] \nabla\cdot\boldsymbol{\sigma}  = &~ \boldsymbol{0}   \hspace{3mm} \rm{in}  \hspace{3mm} \Omega, \label{eqn:linearmoment} \\
-\mathcal{G}_c \ell_o (\nabla)^2 \phi + \left[    \frac{\mathcal{G}_c}{\ell_o}  + 2 H^+ \right] \phi  = &~ 2  H^+ \hspace{3mm} \rm{in}  \hspace{3mm} \Omega, 
\label{eqn:diffeqn}
\end{align}
\label{eqn:strongForm}
\end{subequations}
with the following boundary conditions:
\begin{align}
\left[ (1-\phi)^2 + k_p \right] \boldsymbol{\sigma} \cdot \mathbf{n} =& \overline{\mathbf{t}} \quad \rm{on}  \quad \Gamma_{\rm N}, \nonumber \\
\uu =& \overline{\uu} \quad \rm{on}  \quad \Gamma_{\rm D}, \nonumber \\
\nabla \phi \cdot \mathbf{n}  =& 0 \quad \rm{on}  \quad \Gamma \setminus \Gamma_{\rm c},
\end{align}
where $\bvsig = \frac{\partial \psi(\bveps)}{\partial \bveps}$ is the Cauchy stress tensor and $H^+:= \underset{\tau \in [0, t]}{\rm{max}} \Psi^+(\boldsymbol{\varepsilon}(\boldsymbol{x}, \tau))$. To prevent the crack faces from inter-penetration, \eref{eqn:diffeqn} is supplemented with the following constraint:
\begin{equation}
\forall \boldsymbol{x}: \Psi^+ < \Psi^- \Rightarrow \phi:= 0,
\end{equation}
where,
\begin{equation*}
\Psi^{\pm}(\boldsymbol{\varepsilon}) = \frac{1}{2} \lambda \langle \rm{tr} (\boldsymbol{\varepsilon}) \rangle^{2}_{\pm}  + \mu \rm{tr} (\boldsymbol{\varepsilon}^{2}_{\pm}),
\label{eqn:histpositivenegative}
\end{equation*}
with $\langle \cdot \rangle_{\pm}:= \frac{1}{2}(\pm |\cdot|)$, $\boldsymbol{\varepsilon}_{\pm}:= \sum\limits^3_{I = 1} \langle \varepsilon_I \rangle_{\pm} \mathbf{n}_I \otimes \mathbf{n}_I  $ and $\varepsilon = \sum\limits^{3}_{I = 1} \langle \varepsilon_I \rangle\boldsymbol{n}_I \otimes \mathbf{n}_I,$ where $\varepsilon_I$ and $\mathbf{n}_I$ are the principal strains and the principal strain directions, respectively.

\subsection{Weak form}
Let $\mathscr{W}(\Omega)$ include the linear displacement field and the phase field variable and let ($\mathscr{U}, \mathscr{P})$ and ($\mathscr{V},\mathscr{Q}$), be the trial and the test function spaces:
\begin{subequations}
\begin{align}
(\mathscr{U},\mathscr{V}^0) &=
\left\{ (\uu^h,\vm{v}) \in [ C^0(\Omega)]^d : (\uu,\vm{v}) \in
[ \mathcal{W}(\Omega)]^d \subseteq [ H^{1}(\Omega)]^d \right\},  \\
(\mathscr{P},\mathscr{Q}^0) &= \left\{ (\phi,q) \in
[ C^0(\Omega) ]^d : (\phi, q) \in [ \mathcal{W}(\Omega)]^d \subseteq
[ H^{1}(\Omega) ]^d
\right\}.
\end{align}
\end{subequations}
Let the domain be partitioned into elements $\Omega^h$ and on using shape functions $N$ that span at least the linear space, we substitute the trial and the test functions: $\left\{\vm{u}^h, \vm{\phi}^h\right\} = \sum\limits_I N_I \left\{\vm{u}_I,\phi_I\right\}$ and $\left\{\vm{v}, \theta \right\} = \sum\limits_I N_I \left\{\vm{v}_I,\theta_I\right\}$ into \eref{equation:weakform}, the system of equations can be readily obtained upon applying the standard Bubnov-Galerkin procedure. Find $ \, \uu^h \in \mathscr{U} \, {\rm and} \, \vm{\phi}^h \in \mathscr{P} \, \text{such that, for all} \, \vm{v} \in \mathscr{V}^0 \, {\rm and} \, \theta \in \mathscr{Q}$,
\begin{subequations}
\begin{align}
    \int_{\Omega} \left\{ \left[ (1-\phi)^2+k_p \right]\bm{\sigma}(\vm{u}):\bm{\varepsilon}(\vm{v}) \, \right\} \text{d} \Omega  =& \int_{\Gamma_t} \bm{\bar{t}} \cdot\vm{v}\,\text{d} \Gamma,
\label{equation:elasticity}\\
    \int_{\Omega} \left[ \nabla \theta \, G_c \ell_o\, \nabla \phi + \theta \left( \frac{G_c}{\ell_o} + 2 H^+ \right) \phi \right] \, \text{d} \Omega =& \int_{\Omega} 2 H^+ \theta \, \text{d} \Omega + \int_{\Gamma} \nabla \phi \cdot \vm{n} \, \theta \, \text{d} \Gamma,
    \label{equation:diffusion}
\end{align}
\label{equation:weakform}
\end{subequations}
which leads to the following system of linear equations:
\begin{subequations}
    \begin{align}
    \mathbf{K}^{\rm uu} \vm{u}^h &= \vm{f}^{\rm uu}, \label{eqn:elast_discrete} \\
    \mathbf{K}^{\phi} \vm{\phi}^h &= \vm{f}^{\phi}, \label{eqn:phase_discrete}
\end{align}
\end{subequations}
where
\begin{align*}
    \mathbf{K}^{\rm uu} &= \sum\limits_h \int\limits_{\Omega^h} \Big[(1-\phi)^2+k_p\Big] \mathbf{B}^{\rm T}~\mathbb{D}~\mathbf{B}~\mathrm{d}\Omega, \\
    \mathbf{K}^{\phi} &= \sum\limits_h \int\limits_{\Omega^h} \Bigg[ \mathbf{B}^{\rm T}_{\phi}~G_c \ell_o~\mathbf{B}_{\phi} + \mathbf{N}^{\rm T}\left( \frac{G_c}{\ell_o} + 2 H^+ \right) \mathbf{N} \Bigg]~\mathrm{d}\Omega, \\
    \vm{f}^{\rm uu} &= \sum\limits_h \int\limits_{\Omega^h}  \mathbf{N}^{\rm T} \bm{\bar{t}}~\mathrm{d}\Omega, \\
    \vm{f}^{\phi} &= \sum\limits_h \int\limits_{\Omega^h} \mathbf{N}^{\rm T}~2 H^+~\mathrm{d}\Omega,
\end{align*}
where $\mathbb{D}$ is the material constitutive matrix, $\mathbf{B}= \boldsymbol{\nabla} \mathbf{N}$ is the strain-displacement matrix and $\mathbf{B}_{\phi} = \nabla \mathbf{N}$ is the scalar gradient of the shape function matrix $\mathbf{N}$. The above system of equations are solved by the staggered approach as shown in Algorithm \ref{Algorithm}.


\section{Recovery based error indicator and quadtree decomposition} \label{Sec:Quadtree Meshing}

\subsection{Recovery based error indicator}
The recovery based error indicator proposed by Bordas and Duflot \cite{Bordas2007,Bordas_errror_ind_XFEFM2008} is employed to assess the error and list the elements for refinement. In this method, the enhanced strain field is computed using the standard nodal solution through the eXtended Moving Least Square (XMLS) derivative recovery process. The Moving Least Square (MLS) shape functions of the ${n_{x}}$ points whose domain of influence contain points $\mathbf{x}$ is given by:
\begin{equation}
\label{eqn:MLS shapefunction}
{\mathbf{\Psi}^{\rm{T}}(\mathbf{x})}={\mathbf{p}^{\rm{T}}(\mathbf{x})}{\mathbf{A}^{-1}(\mathbf{x})}{\mathbf{B}(\mathbf{x})}
\end{equation}
where ${\mathbf{p}(\mathbf{x})}= \begin{bmatrix} 1 & x & y \end{bmatrix}$ denotes the reproducing polynomial and
\begin{align*}
{\mathbf{A}(\mathbf{x})} &=\sum_{I=1}^{n_{x}}{w_I(\mathbf{x})}{\mathbf{p}(\mathbf{x}_{I})}{\mathbf{p}^{\rm{T}}(\mathbf{x}_{I})} \nonumber \\
{\mathbf{B}(\mathbf{x})}&=
\begin{bmatrix}
{w_1(\mathbf{x})}{\mathbf{p}(\mathbf{x}_{1})} &
{w_2(\mathbf{x})}{\mathbf{p}(\mathbf{x}_{2})} & ... &  {w_{n_{x}}(\mathbf{x})}{\mathbf{p}(\mathbf{x}_{n_{x}})}
\end{bmatrix}.
\end{align*}
Here, the weight function of a node is calculated by the diffraction method with a circular domain of influence for node $\mathbf{x}_k$. In this work, a fourth order spline is taken as the weighting function:
\begin{equation}
w_k(\mathbf{x}) = 
\begin{cases}
    1- 6{s^2}+8{s^3}-3{s^4} & \text{if} \hspace{0.3cm} \vert s \vert \leq 1  \\ 
	 0\hspace{1cm}               & \text{if} \hspace{0.3cm}\vert s \vert > 1
\end{cases}.
\end{equation}
where $s=\frac{\Vert \mathbf{x}-\mathbf{x}_{k}\Vert}{d_{k}}$ and $d_k$ denotes the support domain of node $\mathbf{x}_k$.  
The enhanced derivative of the shape functions are computed by finding the derivatives of the MLS shape functions (see \eref{eqn:MLS shapefunction}). Using the enhanced derivatives, the enhanced derivatives of the displacements and the enhanced small strain $\mathbf{\varepsilon}_s$ can be found. The error between the enhanced strain field and the standard compatible strain field is considered as the error. The tolerance is considered based on the bulk error criteria, where the fixed fraction of the total error creating elements are refined in the next level. The elements selected for the refinement are chosen by the descending order of their individual elemental error. 

\subsection{Quadtree decomposition}
\begin{figure}[H]
	\centering 
  \subfloat[]{\includegraphics[width=0.4\textwidth]{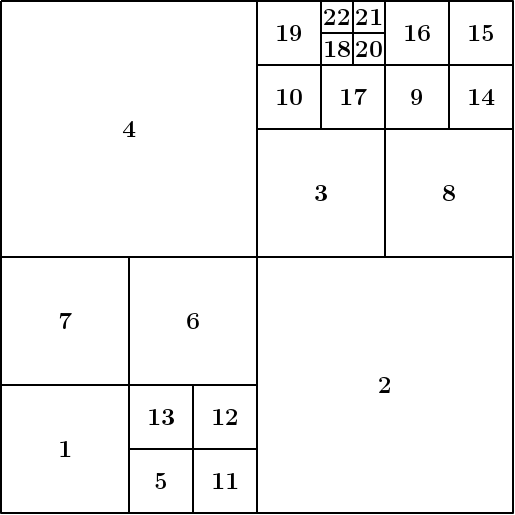}} \hspace{5mm}
  \subfloat[]{\includegraphics[width=0.55\textwidth]{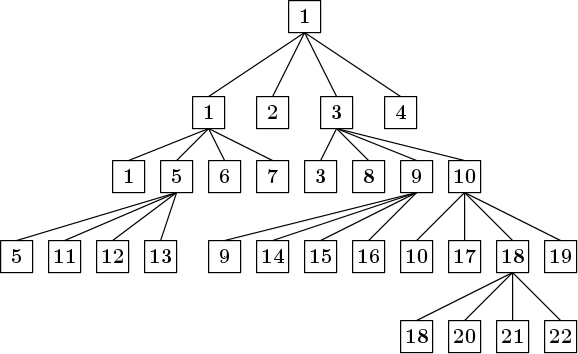}}
\caption{Quadtree decomposition: (a) representative quadtree mesh and (b) tree structure employed to store the mesh details.}
\label{fig:Quadtree mesh details}
\end{figure}

In this meshing technique, one specific criterion named as the stopping criterion is chosen to decide which element requires to be further refined. If the given element does not satisfy the stopping criterion within the user specified tolerance limit, it will be divided into four child elements. A quadtree structure and a mesh with element numbering are shown in \fref{fig:Quadtree mesh details}. This criterion can be a geometry based factor or any error indicator.
\begin{figure}[H]
\centering
\subfloat[]{\includegraphics[scale =0.4]
{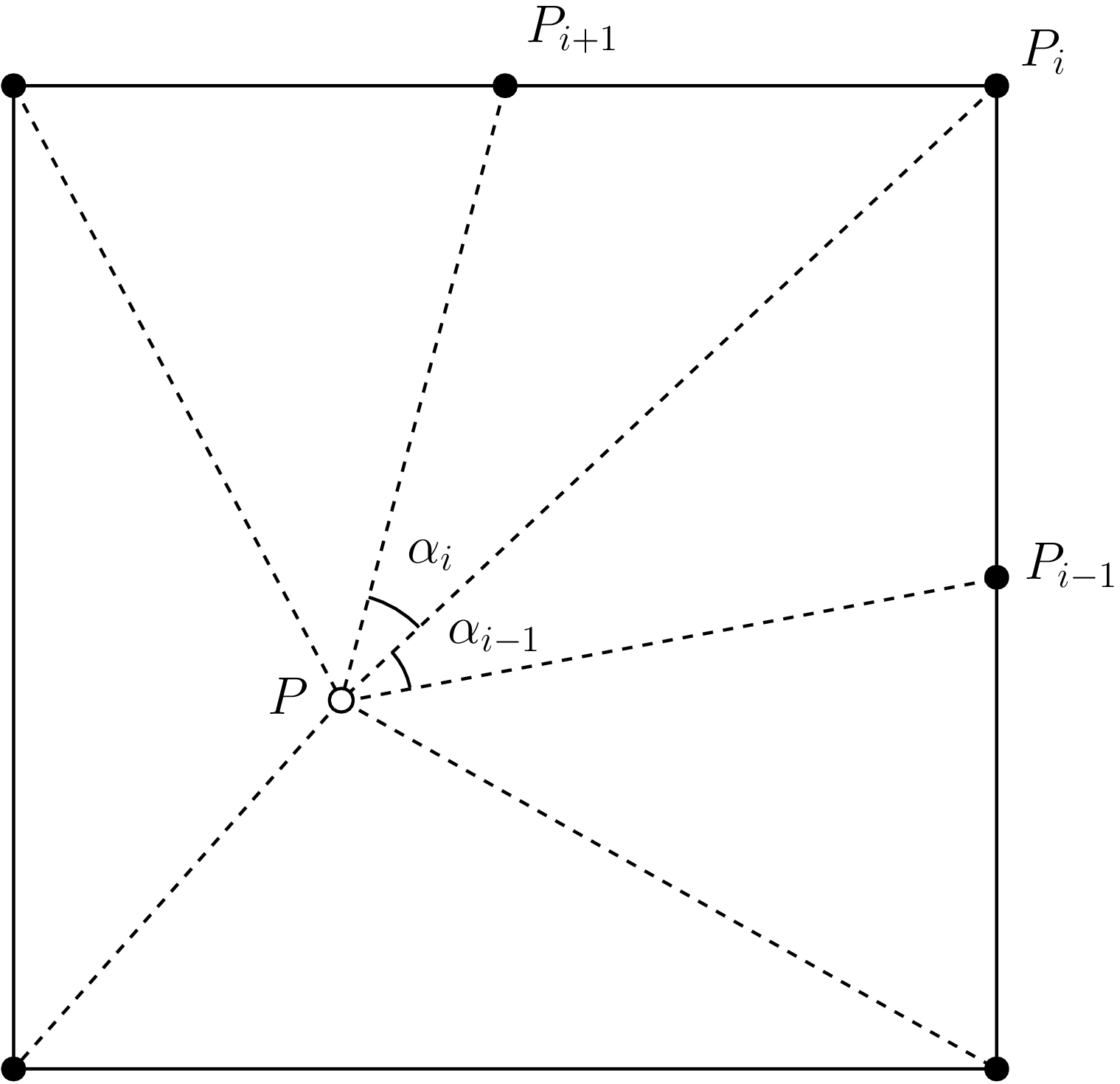}} \hspace{10mm}
\subfloat[]{\includegraphics[scale = 0.6]{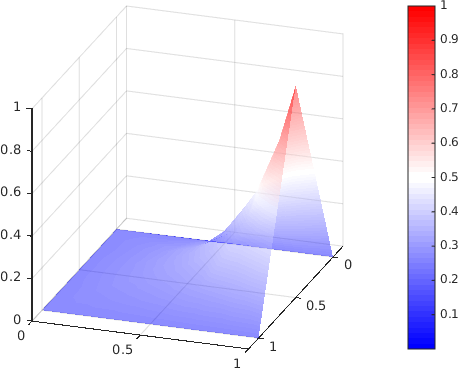}}
\caption {Schematic representation of an element with hanging node and the construction of mean value shape function.} 
\label{fig:Meanvalueshapefunction}
\end{figure}
The above decomposition, leads to elements with hanging nodes and the conventional FE approach cannot handle such elements without additional work. This is because of lack of compatibility. In this work, the elements with hanging nodes are considered as  polygons (see \fref{fig:Meanvalueshapefunction}) and the mean-value coordinates~\cite{Floater2003_meanvalue} are used to approximate the unknown fields. The mean-value coordinates for a point $P(\mathbf{x})$ in an arbitrary polygon is given by:
\begin{align}
N_{i}(\mathbf{x})=\dfrac{\omega_{i}(\mathbf{x})}{\sum_{j=1}^{n} \omega_{i}(\mathbf{x})}, \quad i=1,\cdots,n \nonumber \\
\omega_{i}(\mathbf{x})=\dfrac{\tan(\alpha_{i-1}/2) + \tan(\alpha_{i}/2)}{\Vert \mathbf{x}-\mathbf{x}_{i}\Vert}
\label{meanvalueshapefunction}
\end{align}
where $n$ is the number of nodes in an element, $\mathbf{x}_{i}$ is the coordinate of point ${P}_{i}$ and $\alpha_i$'s are the internal angle. 
\subsection{Solution algorithm}

The staggered solution scheme is used to solve the coupled systems arising from \erefs{equation:weakform}. In the staggered approach, the unknown variables (phase-field variable, $\phi$ and the displacement, $\bm{u}$) are split algorithmically and \erefs{equation:elasticity}-(\ref{equation:diffusion}) are solved in each time step till the convergence as shown in Algorithm \ref{Algorithm}. For a step ($i$), the variables, viz., displacement field ($\bm{u}$), phase-field variable ($\phi$), history variable ($H^+$) are initialized and a converged mesh using the MLS post-errori error estimator in combination with quadtree decomposition is obtained. For the applied displacement ($\Tilde{u}$) at step ($i+1$), the phase-field variable ($\phi^h_{i+1}$) is computed using the \eref{equation:diffusion}. 
Now the displacement field ($\vm{u}^h_{i+1}$) using \eref{equation:elasticity} is computed keeping phase-field variable ($\phi^h_{i+1}$) constant. Then $\psi^+$ $\&$ $\psi^-$ are computed and $H^+$ is updated. The convergence for the displacement and the phase-field variable at the current step and previous step is checked using max \{ $\frac{||\mathbf{u}^h_{i+1}-\mathbf{u}^h_i||}{||\mathbf{u}^h_{i+1}||}$, $\frac{||\phi_{i+1}-\phi_i||}{||\phi_{i+1}|| }$ \} $\leq$ tolerance. 
Once the error in solutions is within the user specific tolerance, 
the mesh and the history variable ($H^+$) are updated. The analysis then proceeds to the next load step.

\begin{algorithm}[H]
Initialize at step $(i)$: 
Converged quadtree mesh, $\vm{u}^h_i$, $\phi^h_i$ and $H^+_i$, $\Tilde{u} = \Delta u$ \\
\For {$\Tilde{u}$=$\Delta u, 2\Delta u, .... , n_{\rm{total}}\Delta u$}{
\While {$|| \vm{u}^h_{i+1} - \vm{u}^h_{i} ||/ ||  \vm{u}^h_{i} ||$ and $|| \phi^h_{i+1} - \phi^h_{i} ||/ || \phi^h_{i} ||$ $\ge$ \rm{tolerance}}{Compute phase-field variable ($\phi^h_{i+1}$) from the \eref{eqn:phase_discrete} \\
Compute displacement field ($\vm{u}^h_{i+1}$) from \eref{eqn:elast_discrete}\\ 
Compute $\psi^+$ and $\psi^-$ and update $H^+$ \\
\If{Iteration $\ge$ numIter}{Update mesh using post-processing error estimator and quadtree decomposition.}
} 
Update history variable ($H^+$) \\
Update mesh using post-processing error estimator and quadtree decomposition.
} 
\caption{{\bf staggered solution algorithm for adaptive PFM} \label{Algorithm}}

\end{algorithm}

\section{Numerical examples}
\label{Section:numericalexamples}
In this section, we test the proposed adaptive PFM on the standard PFM literature numerical examples. In order to show the robustness and accuracy of the proposed method, we compare the number of elements required to represent the phase-field variable ($\phi$) in the standard PFM literature. We start with a square domain containing a straight edge crack, which is well studied in the PFM literature, see e.g. \cite{Bourdin2008,Miehe2010,Ambati2015ductile,Ambati2015}. Section \ref{section:uniaxialtension} and \ref{section:shear} presents a plane strain domain with an edge crack specimen subjected to tension and shear, respectively. Section \ref{section:lshape} presents and validates the mix-mode crack propagation in L-shape panel. The length scale parameter $\ell_o$ is assumed to be 2$h$ if  not stated otherwise, where $h$ is the minimum size of the element. The numerical stability parameter $k_p$ is assumed to be 1 $\times 10^{-6}$ in all the numerical examples. The total degrees of freedom (Dof) in all the examples are represented as, ndofp $\times$ numnode + ndofu $\times$ numnode, where numnode is the number of nodes, ndofp = 1 and ndofu = 2 are the degrees of freedom per node for the phase-field variable ($\phi$) and the displacement, respectively.    

\subsection{Uniaxial tension (Mode-I fracture) test}
\label{section:uniaxialtension}

\begin{figure}[H]
    \centering
    \subfloat[]{\includegraphics[scale = 0.5]{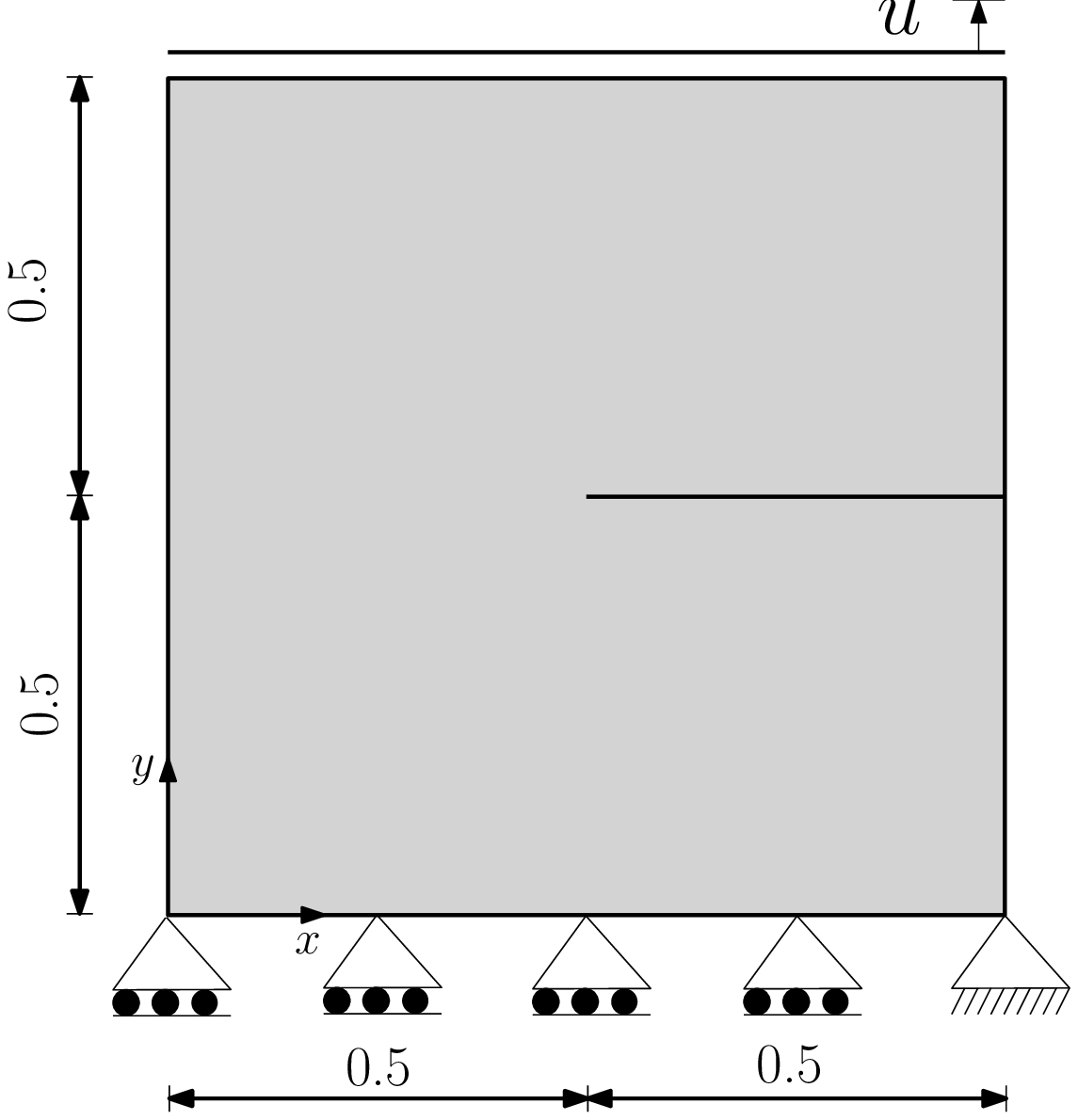}\label{fig:tension}} \hspace{20mm}
    \subfloat[]{\includegraphics[scale = 0.5]{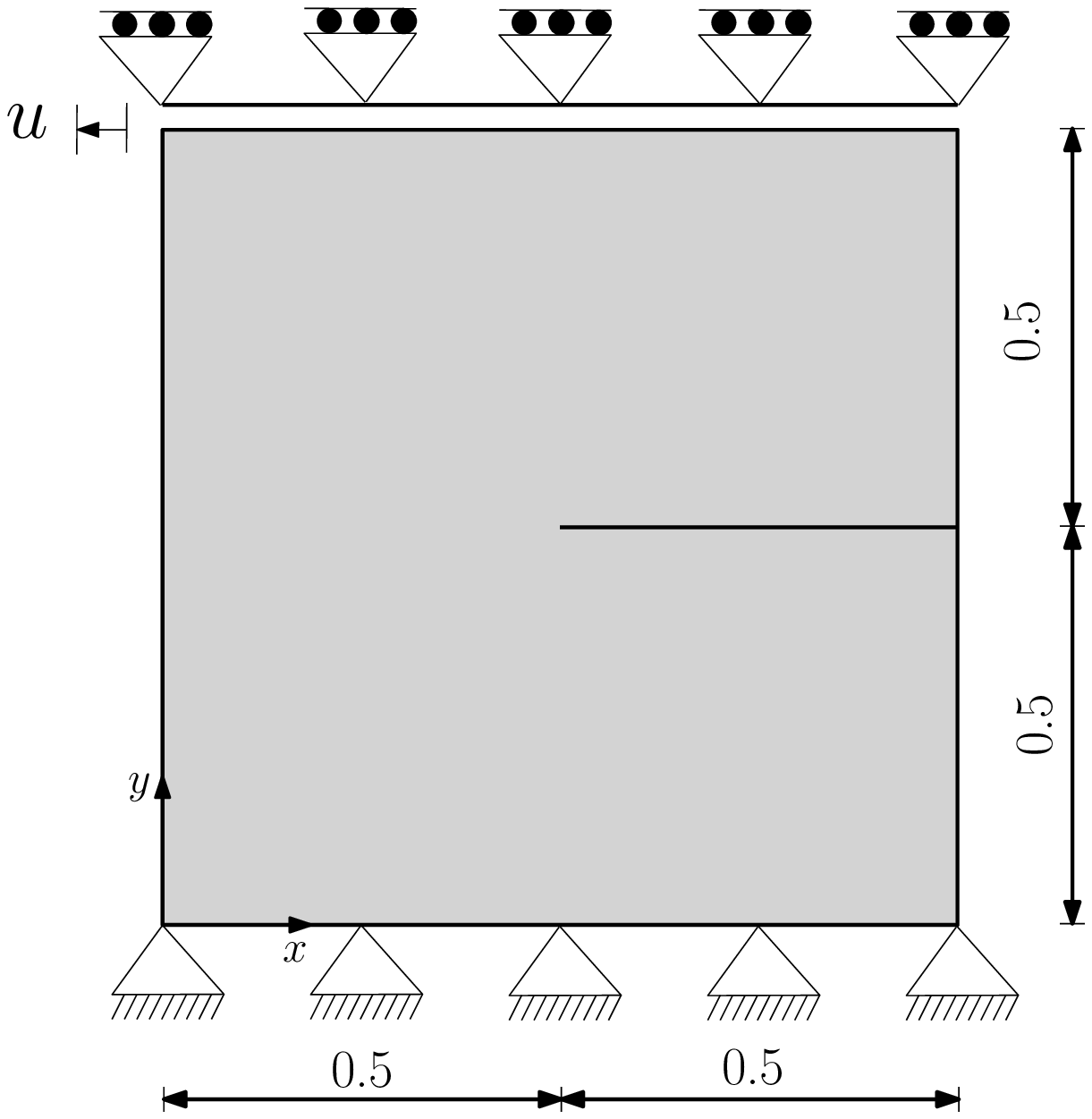}\label{fig:shear}}
    \caption{A plate with an edge crack subjected to (a) tension and (b) shear.}
    \label{fig:domain_edgecrack}
\end{figure}

As shown in \fref{fig:domain_edgecrack}, we consider a plane strain domain with a geometrically induced edge crack. The displacement is applied to the top edge in $y$-direction and the bottom edge is fixed (see \fref{fig:tension}) in order to simulate the mode-I fracture. Displacement is applied with an increment of $\Delta u$ = 1 $\times 10^{-5}$ mm up to $u=$ 5$\times$10$^{-3}$ mm and $\Delta u=$ 1$\times$10$^{-6}$ mm up to failure of the specimen. Similar to the work of \citet{Miehe2010}, the material properties of the specimen are chosen as $\lambda=$ 121.15 kN/mm$^2$, $\mu=$ 80.77 kN/mm$^2$, $G_c=$ 2.7 MPa mm, where $\lambda$ and $\mu$ are Lam$\acute{e}$ constants and $G_c$ is the critical fracture toughness. 

\begin{figure}[H]
    \centering
    \subfloat[]{\includegraphics[scale =0.8]{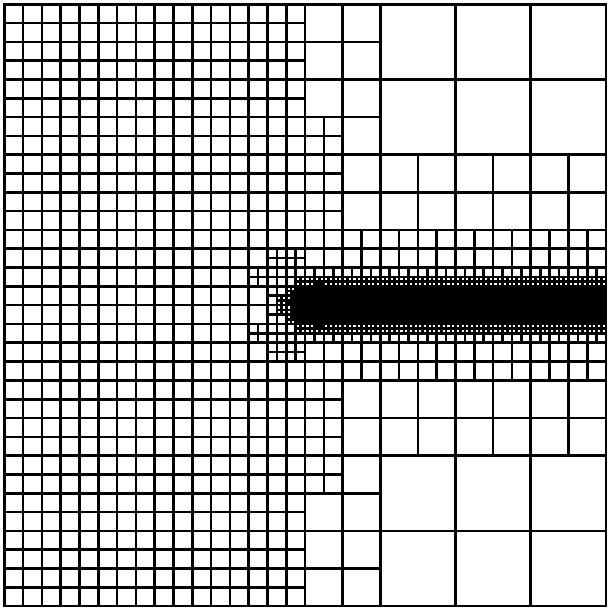}} \hspace{5mm}
    \subfloat[]{\includegraphics[scale =0.8]{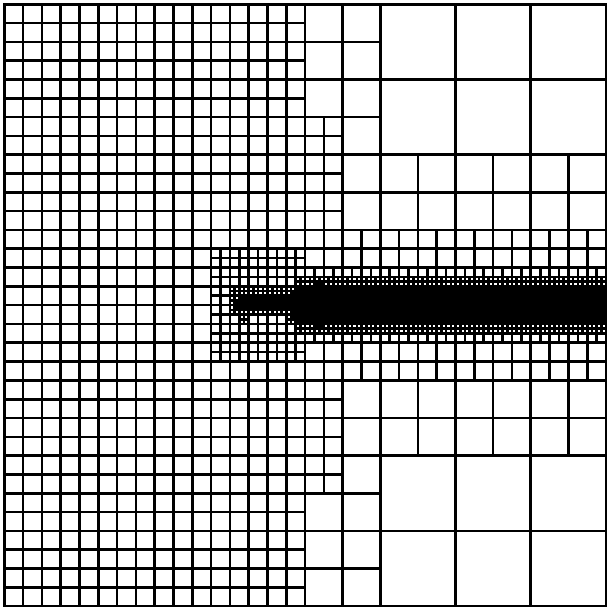}} \hspace{5mm}
    \subfloat[]{\includegraphics[scale =0.8]{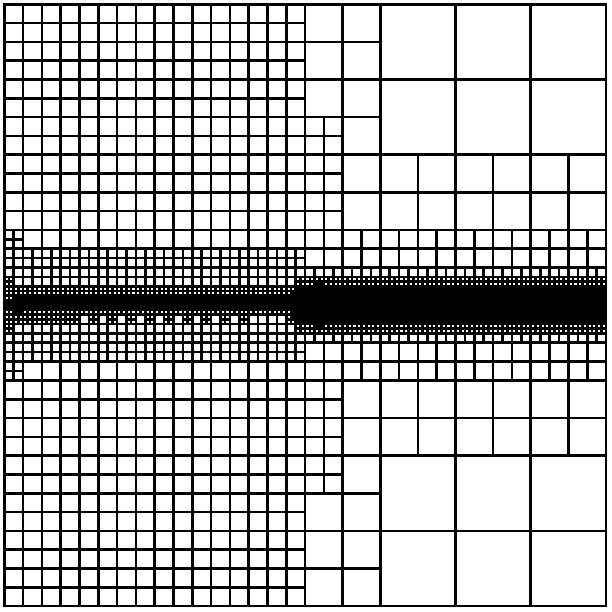}}

    \subfloat[]{\includegraphics[scale =0.19]{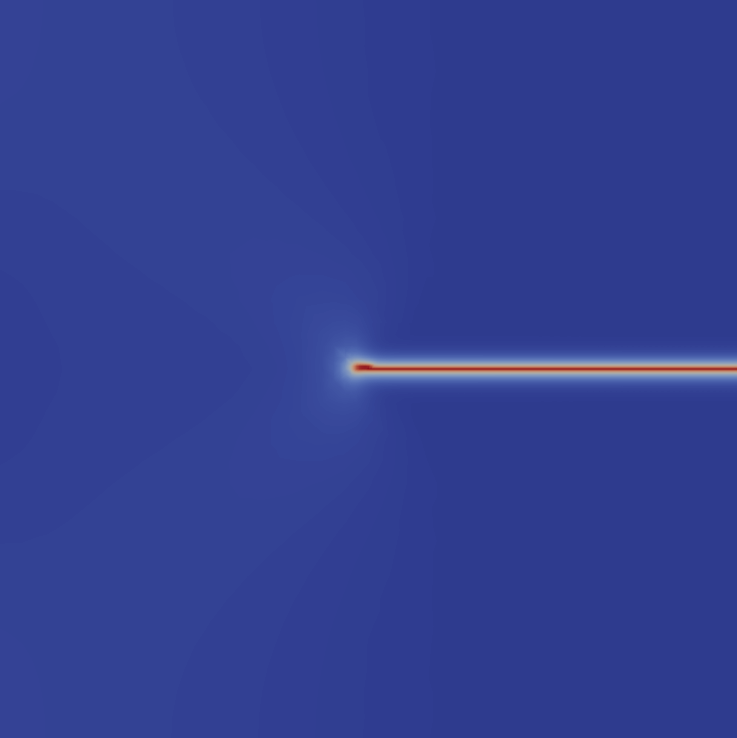}} \hspace{5mm}
    \subfloat[]{\includegraphics[scale =0.19]{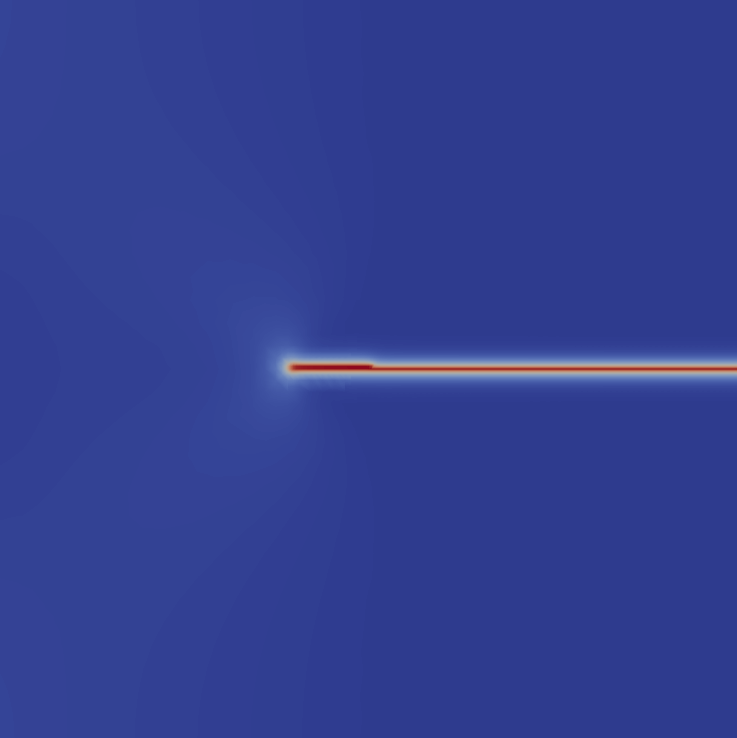}} \hspace{5mm}
    \subfloat[]{\includegraphics[scale =0.19]{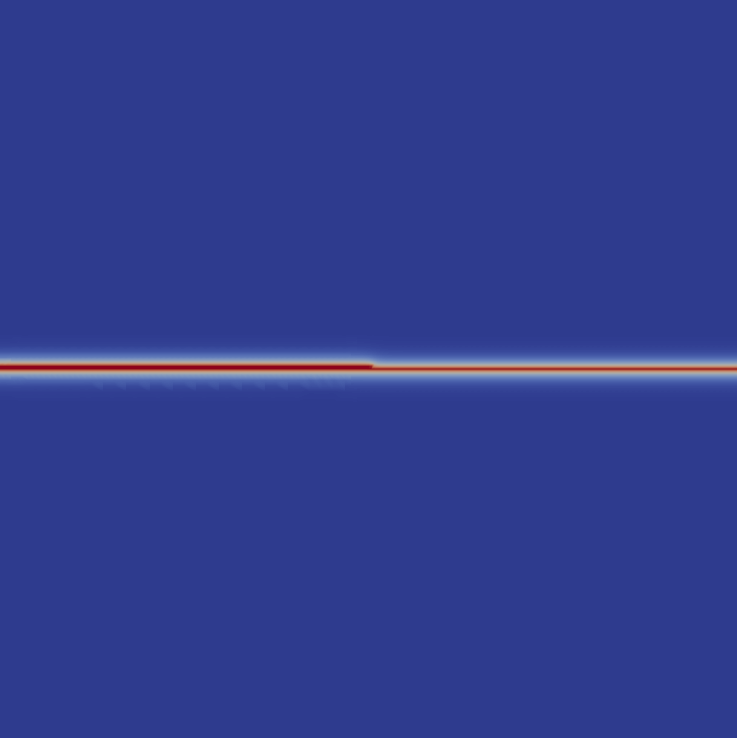}}
    \caption{A plate with an edge crack subjected to tension: crack trajectory at remote displacement (a) 0.0057 mm (elements, 3,181 and Dofs, 10,557) (b) 0.0058 mm (elements, 3,382 and Dofs, 11,214) (c) 0.0059 mm (elements, 4,300 and Dofs, 14,262). }
    \includegraphics[scale = 0.28]{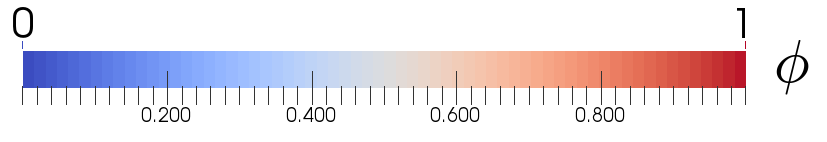}
    \label{fig:tension_phi}
\end{figure}

The initial converged discretization of the domain is obtained using the quadtree decomposition and the MLS posterrori error estimator as explained in the Section \ref{Sec:Quadtree Meshing}. The converged discretization consists of 3,160 quadtree elements and 10,482 Dofs. The crack starts propagating once the stress at the crack tip reaches critical stress. \fref{fig:tension_phi} shows the snapshots of the crack propagation and the corresponding discretization with the applied displacement. The density of the discretization concentrates around the vicinity of the crack tip. The combination of quadtree decomposition and posterrori error estimator strategy leads to fine discretization locally and coarse globally. Moreover, the adaptive strategy automatically adapts the advancement of the crack propagation. A total of 4,300 elements and 14,262 Dofs are required at the end of the simulation. This feature of the proposed adaptive PFM drastically reduced the required elements and Dofs to represent the phase-field variable ($\phi$).      

\begin{figure}[H]
    \centering
    \subfloat[]{\includegraphics[scale = 0.45]{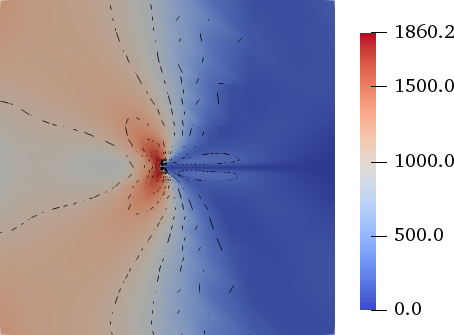}\label{fig:stress_tensiona}} \hspace{7mm}
    \subfloat[]{\includegraphics[scale = 0.36]{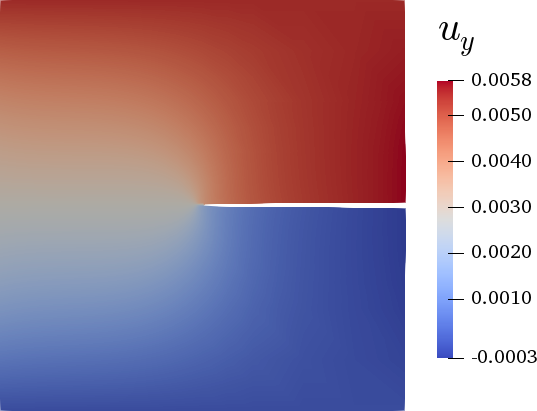}\label{fig:stress_tensionb}}
    
    \subfloat[]{\includegraphics[scale = 0.45]{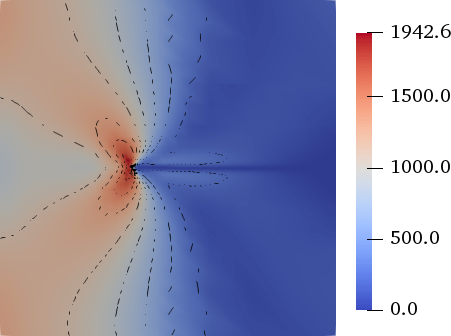}\label{fig:stress_tensionc}} \hspace{7mm}
    \subfloat[]{\includegraphics[scale = 0.36]{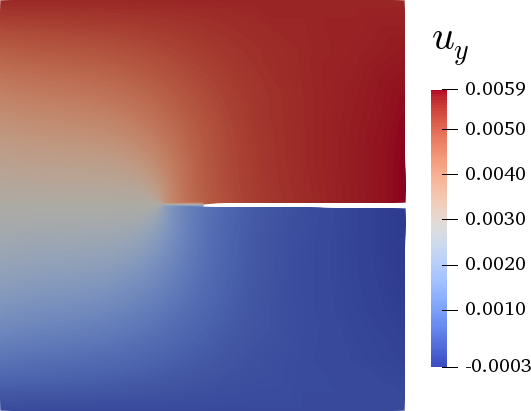}\label{fig:stress_tensiond}}
    
    \subfloat[]{\includegraphics[scale = 0.45]{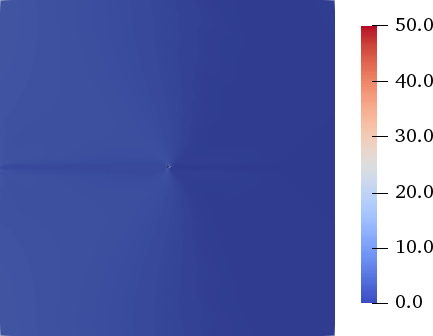}\label{fig:stress_tensione}} \hspace{7mm}
    \subfloat[]{\includegraphics[scale = 0.36]{./Figures/Tension_deformation_650.png}\label{fig:stress_tensionf}}
    \caption{A plate with an edge crack subjected to tension: von-Mises stress distribution and deformed configuration at remote displacement (a, b) 0.0057 mm (c, d) 0.0058 mm (e, f) 0.0059 mm [The deformation is scaled by a factor 2 for visualization].}
    \label{fig:stress_tension}
\end{figure}

\begin{figure}[H]
    \centering
    \includegraphics[scale = 0.8]{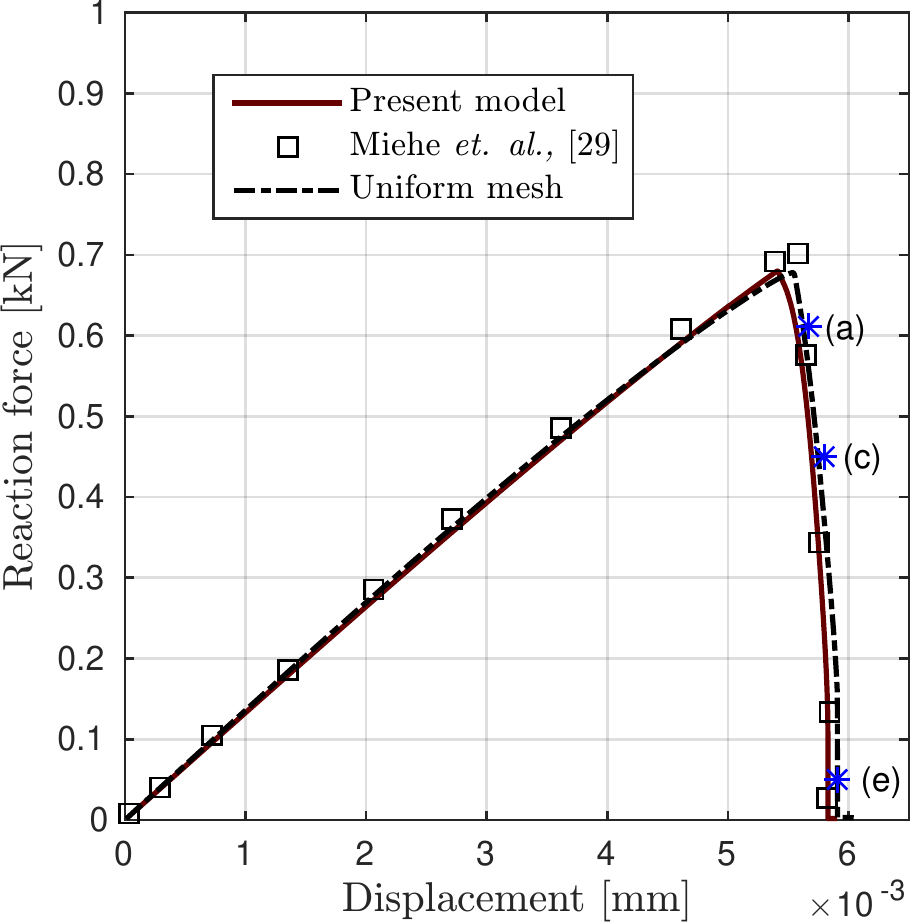}
    \caption{Edge crack plate subjected to tension: load displacement curve  [Points A, C and E are the locations of phase field parameter, stress and displacement plots].}
    \label{fig:load_disp_tension}
\end{figure}

Apart from the phase-field variable ($\phi$) map,  the crack can be also observed in the stress map and the displacement field. The stresses are localized around the crack tip and displacement shows discontinuity around the crack. \fref{fig:stress_tension} shows the von-Mises stress distribution and the corresponding deformed configuration at the locations (A, B and C) marked in the \fref{fig:load_disp_tension}. The von-Mises stress concentrates around the moving crack tip as shown in \frefs{fig:stress_tensiona} and \ref{fig:stress_tensionc}. The displacement in the $y$-direction shows discontinuity, see \frefs{fig:stress_tensionb}, \ref{fig:stress_tensiond} and \ref{fig:stress_tensionf}.

\fref{fig:load_disp_tension} compares the load-displacement response obtained from the present model with the work of \citet{Miehe2010} and the domain with uniform discretization. The uniform dicretization of the domain is obtained using 150 $\times$ 150 divisions,  which results in 22,500 Q4 elements and the minimum element size 0.0066 mm. The characteristic length scale $\ell_o$ is set as 0.0133, which is 2 times minimum size of the element. \citet{Miehe2010} uses 20,000 three noded triangular elements with characteristic length scale $\ell_o$=0.015. The similar results have been obtained in the literature for example \citet{Ambati2015} and \citet{Hirshikesh2018} which uses 12,735 four node quadrilateral elements and 30,546 three node triangular elements, respectively. The presented approach requires lesser number of elements than the uniform discretization and the literature, thus proving the efficacy of the proposed approach. 

\subsection{Shear test (Mode-II fracture)}
\label{section:shear}
In this example we investigate the mode-II fracture. The same specimen used in the previous example, Section \ref{section:uniaxialtension} is used and the incremental displacement $\Delta u = 1 \times 10^{-5}$ mm is applied in the $x$-direction at the top edge till the complete fracture of the specimen while keeping vertical displacement of the top edge fixed, see \fref{fig:shear}. The displacement at the bottom edge is also constrained. The length scale parameter $\ell_o$ is set to 0.01 mm to match the results to the reference solutions \cite{Ambati2015,Hirshikesh2018}. 

\begin{figure}[H]
    \centering
    \subfloat[]{\includegraphics[scale =0.8]{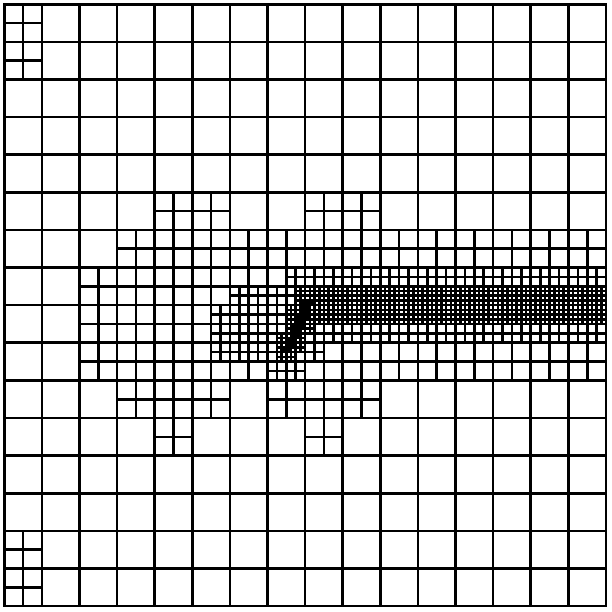}} \hspace{5mm}
    \subfloat[]{\includegraphics[scale =0.8]{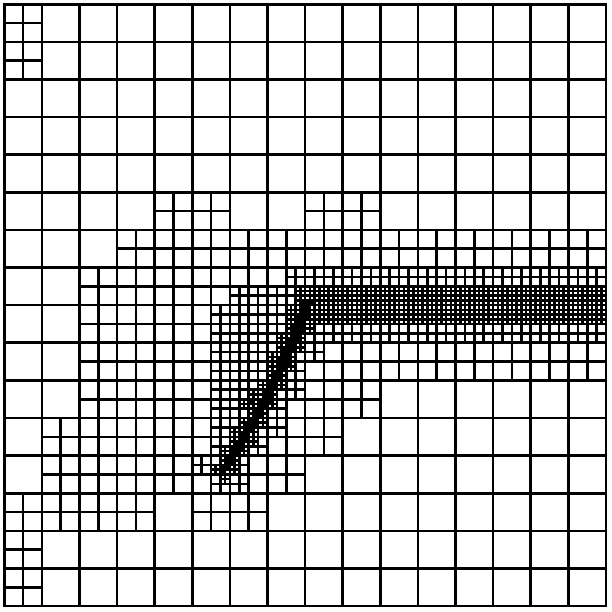}} \hspace{5mm}
    \subfloat[]{\includegraphics[scale =0.8]{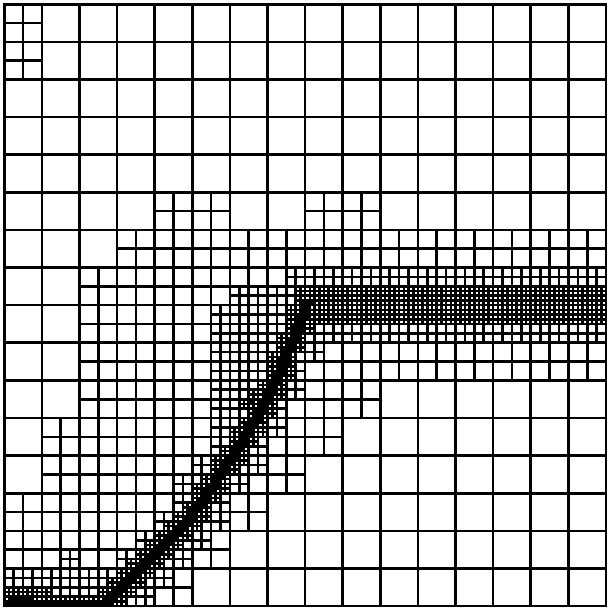}}

    \subfloat[]{\includegraphics[scale =0.19]{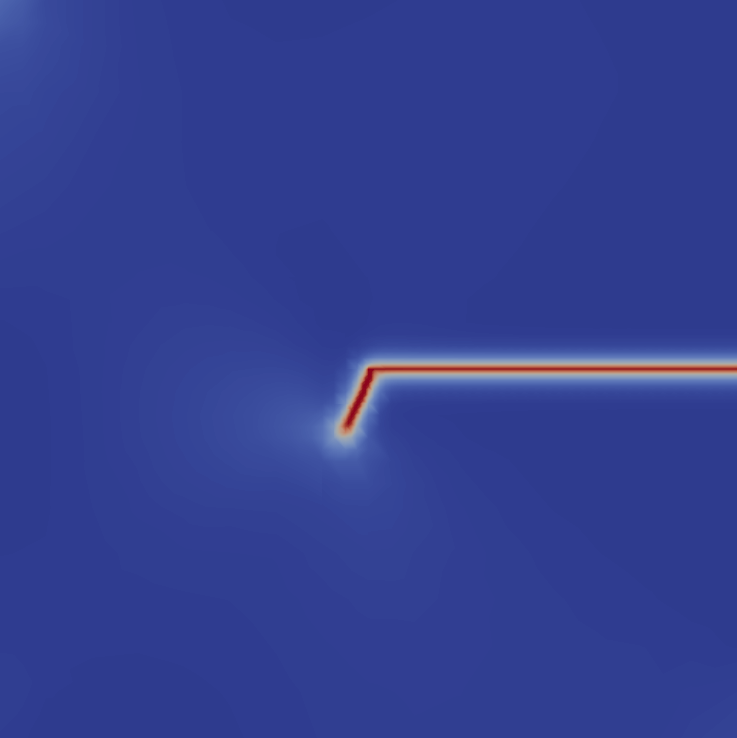}} \hspace{5mm}
    \subfloat[]{\includegraphics[scale =0.19]{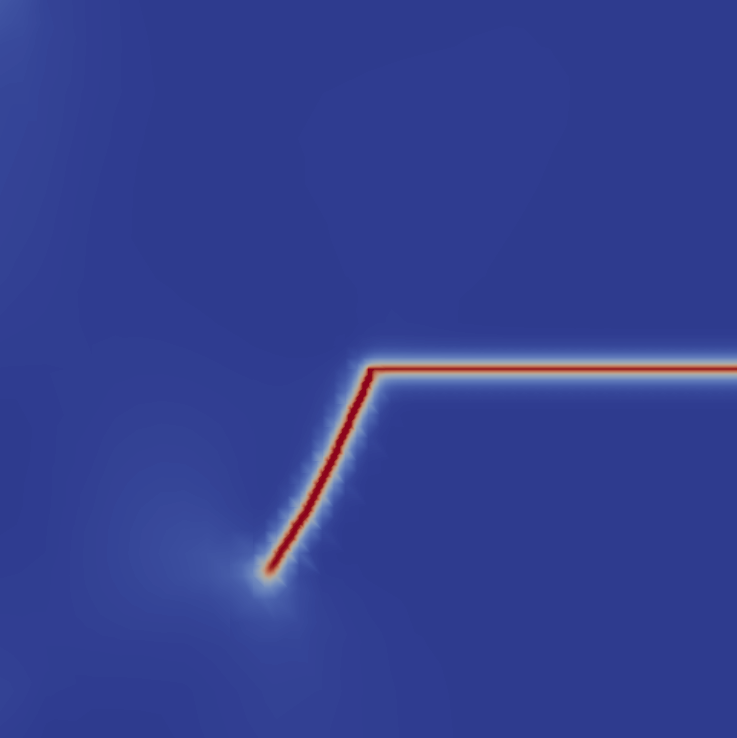}} \hspace{5mm}
    \subfloat[]{\includegraphics[scale =0.19]{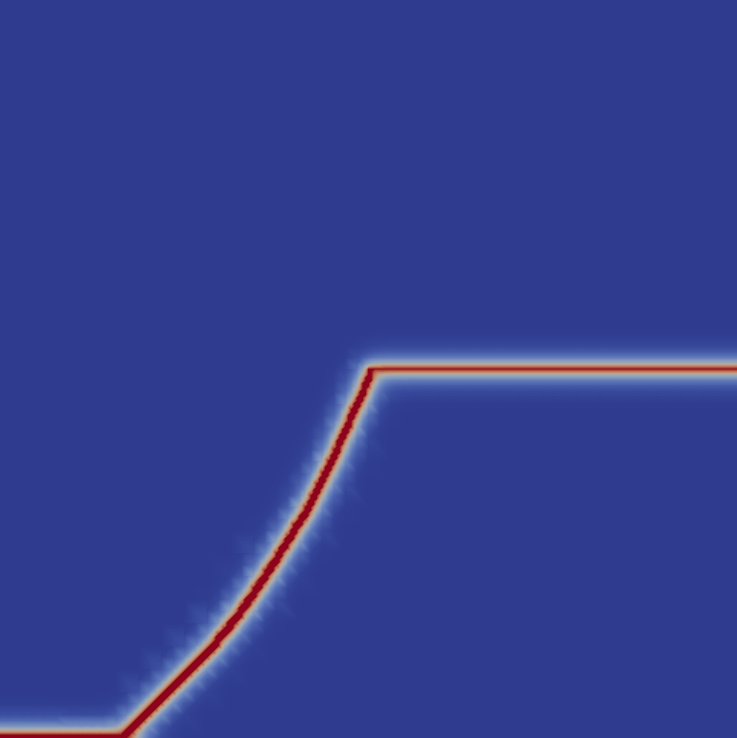}}
    \caption{A plate with an edge crack subjected to shear: the crack trajectory and the corresponding discretization at remote displacement (a) 0.0125 mm (elements, 1,240 and Dofs, 4,371) (b) 0.015 mm (elements, 1,663 elements and Dofs, 5,853) (c) 0.0249 mm (elements, 2,368 elements and Dofs, 8,349).}
    \includegraphics[scale = 0.28]{./Figures/legend_h.png}
    \label{fig:phi_shear}
\end{figure}

\begin{figure}[H]
    \centering
    \includegraphics[scale = 0.8]{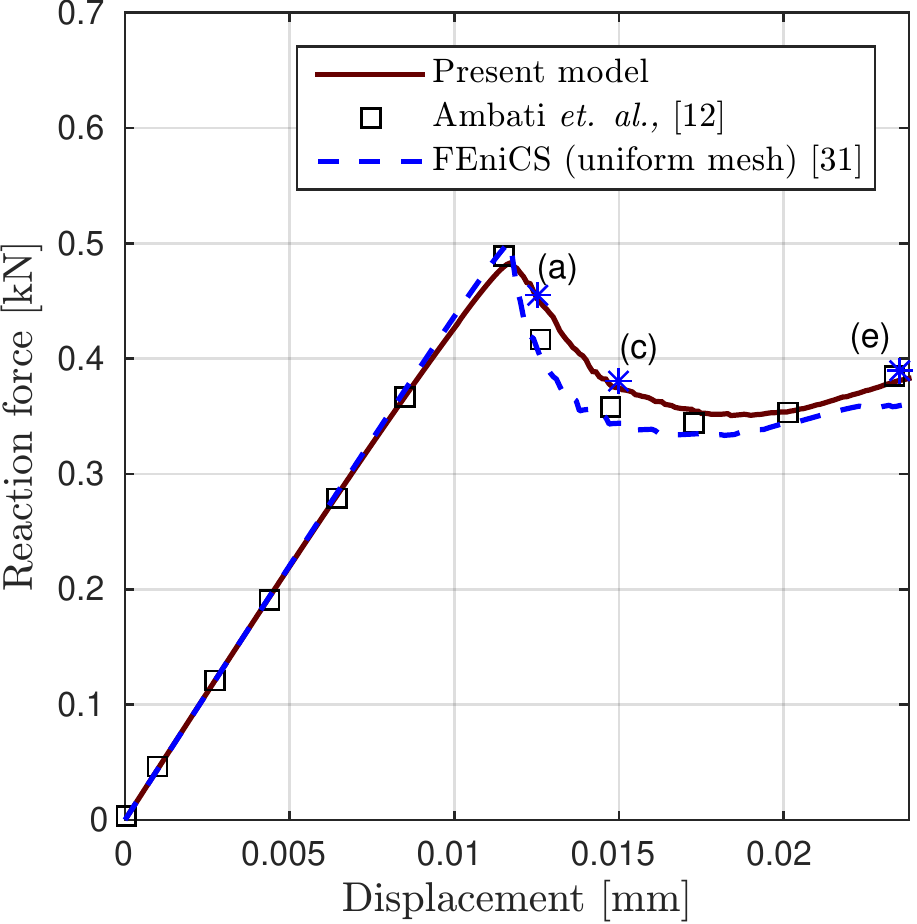}
    \caption{A plate with an edge crack subjected to shear: load displacement curve [Points A, C and E are the locations of phase field parameter, stress and displacement plots].}
    \label{fig:load_disp_shear}
\end{figure}

The simulation starts with the quadtree decomposition which results 1,060 elements and 3,732 Dofs. The crack propagation trajectory and the corresponding dicretization at the locations marked on the \fref{fig:load_disp_shear} are shown in \fref{fig:phi_shear}. The number of elements increases from 1,060 to the 2,368 at the final fracture of the specimen. This is one order small as compared to the global refinement, see for example 30,000 triangular elements \cite{Miehe2010}, 20,592 \cite{Ambati2015} and 40,000 \cite{PATIL2018674} four noded quadrilateral elements. The elements required in the proposed approach are also less than the adaptive PFM 8,664 \cite{PATIL2018254}, local moving XPFM 7,392 \cite{PATIL2018674} and 4,466 \cite{Nagaraja2018}.  

\begin{figure}[H]
    \centering
    \subfloat[]{\includegraphics[scale = 0.35]{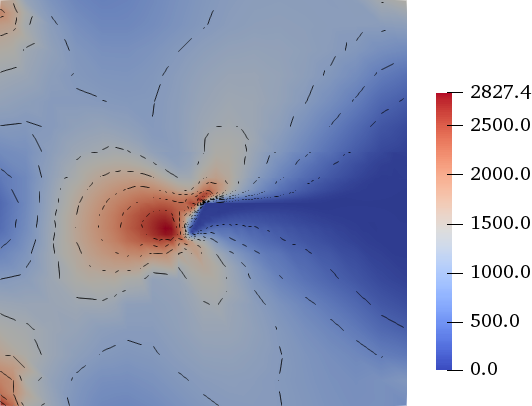}\label{fig:shear_stressa}} \hspace{10mm}
    \subfloat[]{\includegraphics[scale = 0.36]{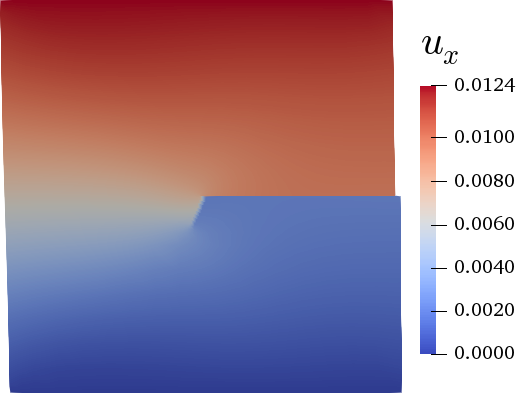}\label{fig:shear_stressb}}
    
    \subfloat[]{\includegraphics[scale = 0.35]{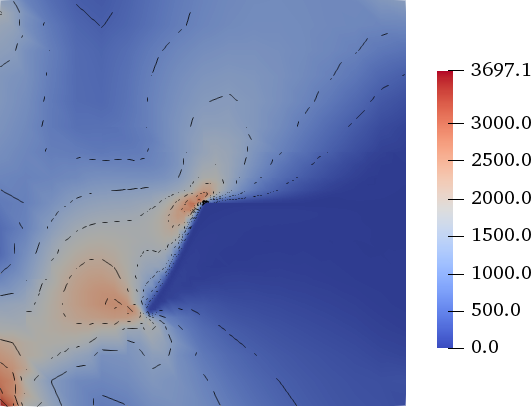}\label{fig:shear_stressc}} \hspace{10mm}
    \subfloat[]{\includegraphics[scale = 0.36]{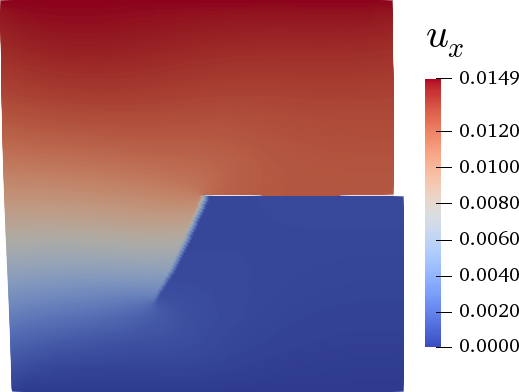}\label{fig:shear_stressd}}
    
    \subfloat[]{\includegraphics[scale = 0.35]{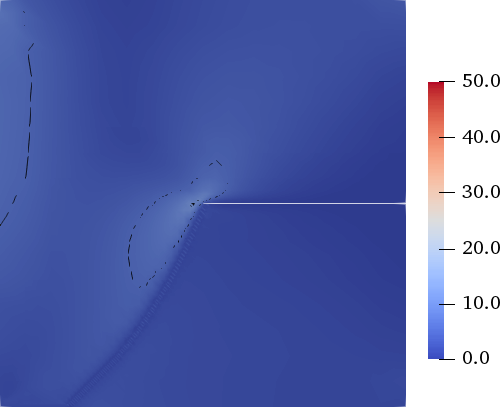}\label{fig:shear_stresse}} \hspace{10mm}
    \subfloat[]{\includegraphics[scale = 0.36]{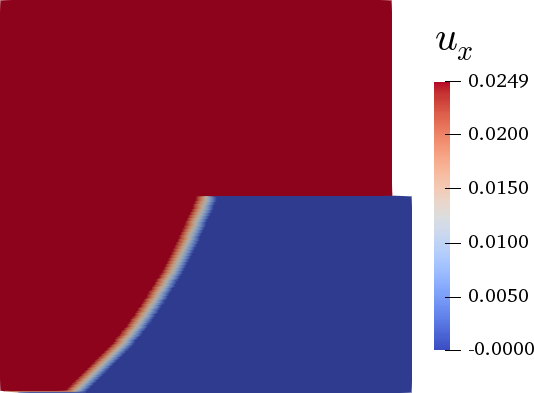}\label{fig:shear_stressf}}
    \caption{A plate with an edge crack subjected to shear: von-Mises stress distribution and deformed configuration at remote displacement (a, b) 0.0124 mm (c, d) 0.0149 mm (e, f) 0.0249 mm [The deformation is scaled by a factor 2 for visualization].}
    \label{fig:stress_shear}
\end{figure}

\fref{fig:stress_shear} shows the von-Mises stress map and the corresponding deformed configuration. The von-Mises stresses are highly localized on one side of the crack tip, which indicates mode-II fracture, see \fref{fig:shear_stressa}. Similar to the tension test, the deformed configuration shows the discontinuity in the displacement field around the crack as shown in \frefs{fig:shear_stressb}, \ref{fig:shear_stressd} and \ref{fig:shear_stressf}. \fref{fig:load_disp_shear} compares the load-displacement with the \citet{Ambati2015} and the \citet{Hirshikesh2018}. The results show very good agreement with the reference solutions. 

\subsection{L-shape panel under mixed-mode fracture}
\label{section:lshape}

\begin{figure}[H]
    \centering
    \subfloat[]{\includegraphics[scale = 0.4]{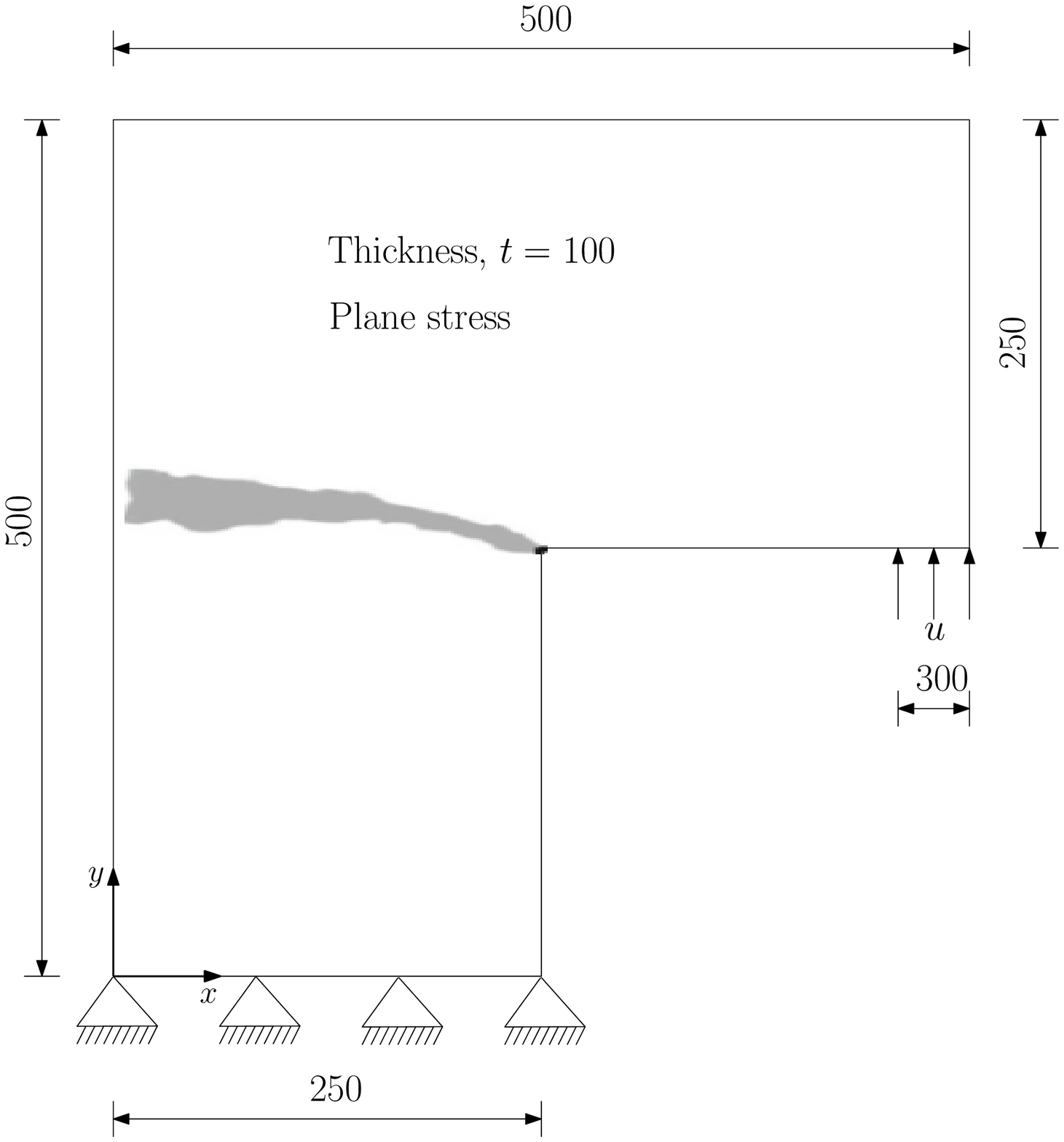}\label{fig:lshape_domain}} \hspace{5mm}
    \subfloat[]{\includegraphics[scale = 0.75]{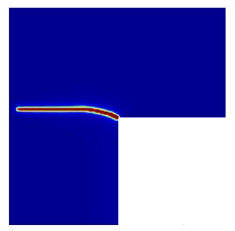}} 
    \caption{L-shape panel specimen: (a) domain, boundary conditions and experimental crack trajectory from \cite{wrinklerB2001} (b) phase field crack propagation path \cite{Ambati2015} [grey region is the experimental crack path and all dimensions are in mm] .}
    \label{fig:lshape_domain_full}
\end{figure}

In the last example, a mixed-mode fracture problem is investigated. We consider a L-shape panel test which is experimentally investigated by \citet{wrinklerB2001}. The geometry and the boundary conditions are depicted in \fref{fig:lshape_domain}. The thickness of the specimen is $t=$ 100 mm. As shown in \fref{fig:lshape_domain} the numerical simulation is displacement controlled and the reaction force is calculated. The incremental displacement $\Delta u = 2.0 \times 10^{-3}$ mm is applied in the $y$-direction at loading point till the fracture of the specimen while both the vertical and the horizontal displacements are constrainted at the bottom edge. Following the work of \citet{wrinklerB2001} the material properties are chosen as Young's modulus $E=25.85$ GPa, Poisson's ratio $\nu = 0.18$ and the critical fracture toughness $G_c = 95$ N/m. The length scale parameter $\ell_o$ is set as 0.2 mm. 

\begin{figure}[H]
    \centering
    \subfloat[]{\includegraphics[scale =0.8]{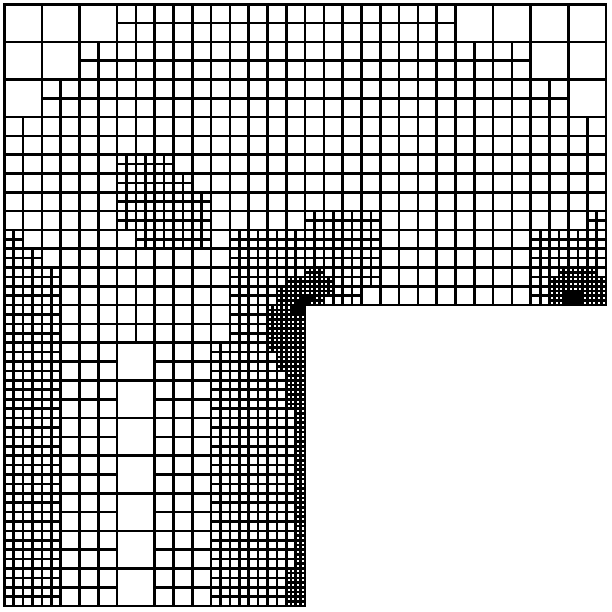}} \hspace{5mm}
    \subfloat[]{\includegraphics[scale =0.8]{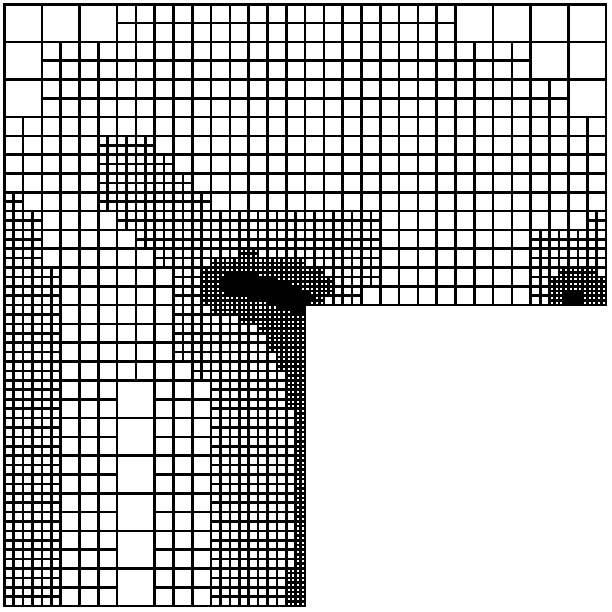}} \hspace{5mm}
    \subfloat[]{\includegraphics[scale =0.8]{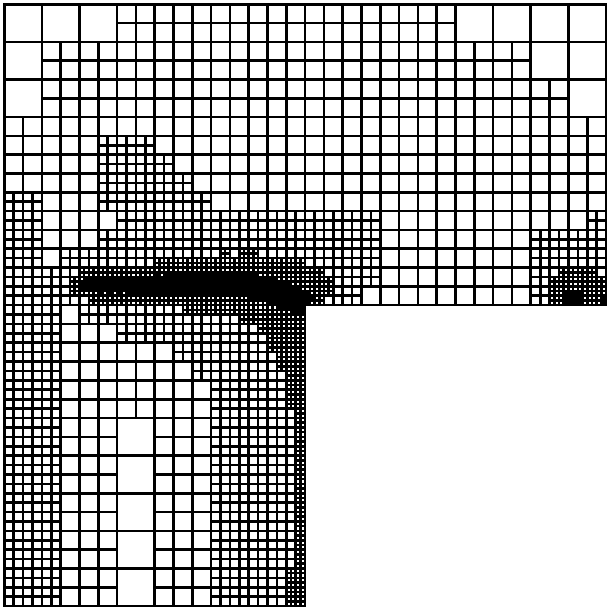}}
    
    \subfloat[]{\includegraphics[scale =0.19]{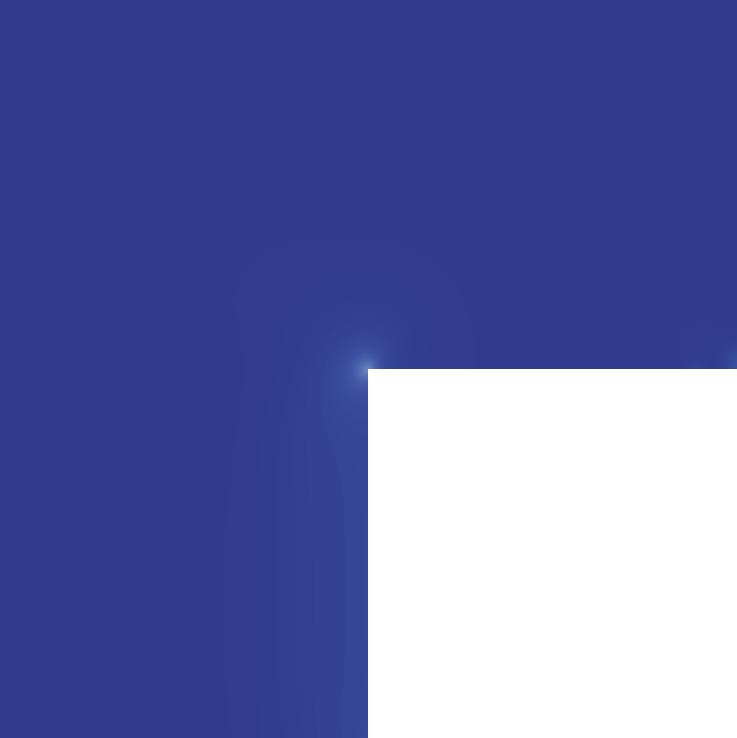}} \hspace{5mm}
    \subfloat[]{\includegraphics[scale =0.19]{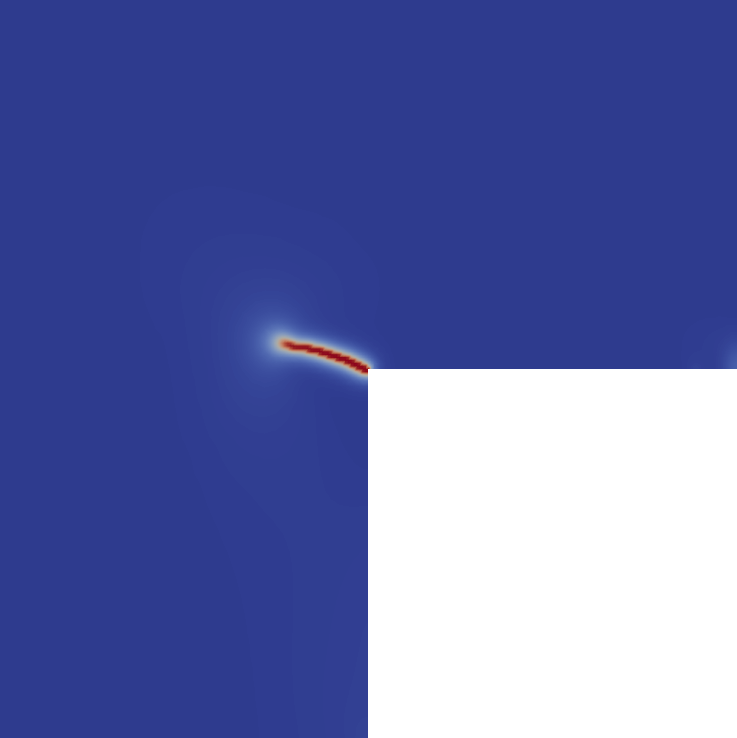}} \hspace{5mm}
    \subfloat[]{\includegraphics[scale =0.19]{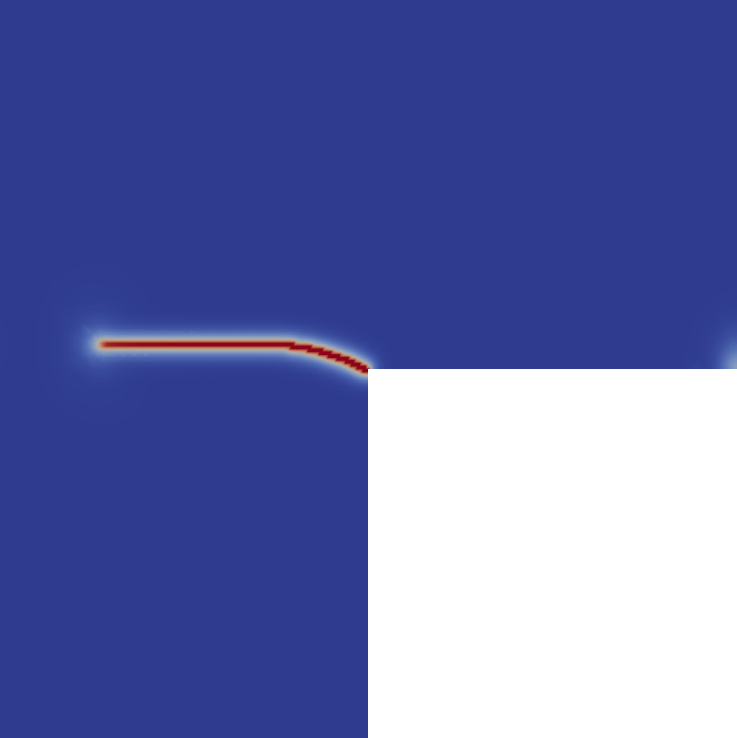}\label{fig:final_lshape}}
    \caption{L-shape panel specimen: crack trajectory and the corresponding discretization at the remote displacement (a) 0.23 mm (elements, 1,701 and Dofs, 5,745)  (b) 0.29 mm (elements, 2,553 and Dofs, 8,439) (c) 0.49 mm (elements, 2,916 and Dofs, 9,636).}
\includegraphics[scale = 0.28]{./Figures/legend_h.png}
    \label{fig:phi_lshape}
\end{figure}

The specimen is discretized initially with 1,617 quadtree elements and Dofs, 5,484. The crack initiates at the corner and starts propagating in curved path. The predicted crack propagation trajectory and the corresponding discretization at different load steps are shown in \fref{fig:phi_lshape}. The discretization automatically adapts as the crack propagates. The number of elements increases from 1,617 to 2,916 at the end of the simulation. The final crack path is compared with the experimental \cite{wrinklerB2001} and the hybrid PFM model \cite{Ambati2015}, see \fref{fig:final_lshape} and \fref{fig:lshape_domain_full}. The simulated crack path is in very good agreement with the experimental results and the hybrid PFM. The number of elements required in the proposed approach is less than the standard PFM 9,650 quadrilateral elements with local mesh refinement \cite{Ambati2015} and the adaptive PFM, 16,524 \cite{PATIL2018254} and 14,792 \cite{PATIL2018674}.

\begin{figure}[H]
    \centering
    \subfloat[]{\includegraphics[scale = 0.24]{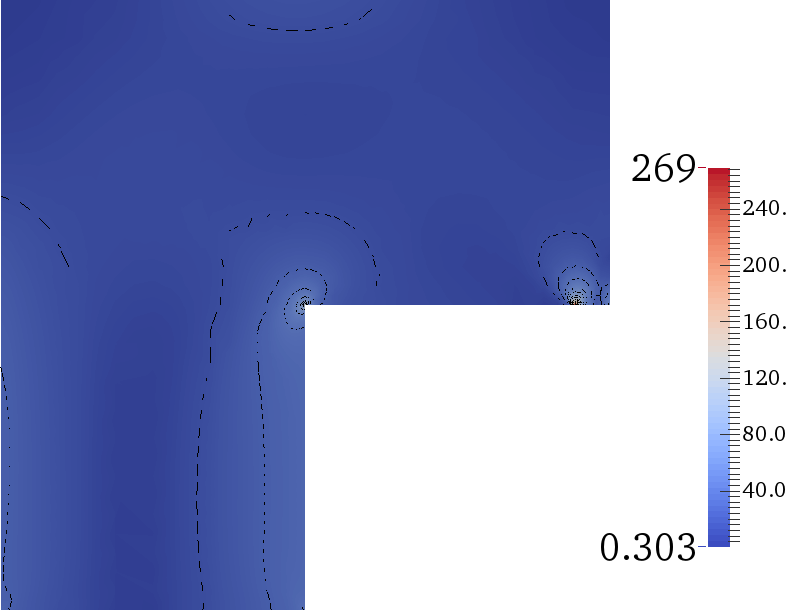}\label{fig:lshape_stressa}} \hspace{10mm}
    \subfloat[]{\includegraphics[scale = 0.35]{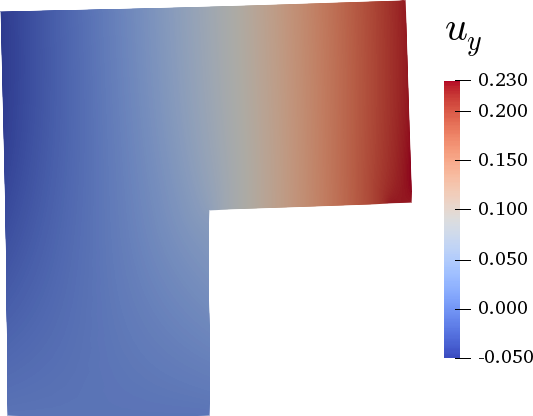}\label{fig:lshape_stressb}}
    
    \subfloat[]{\includegraphics[scale = 0.24]{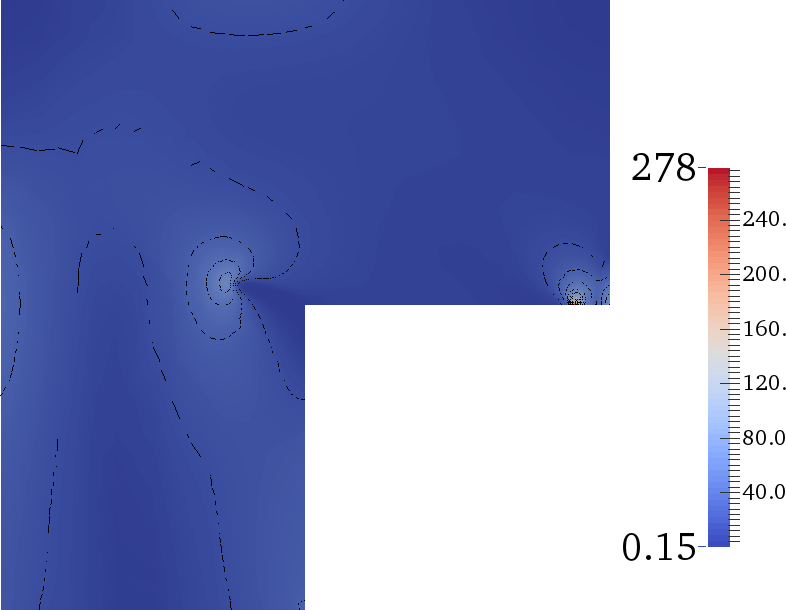}\label{fig:lshape_stressc}} \hspace{10mm}
    \subfloat[]{\includegraphics[scale = 0.4]{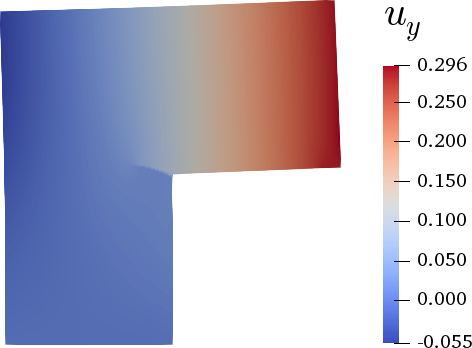}\label{fig:lshape_stressd}}
    
    \subfloat[]{\includegraphics[scale = 0.24]{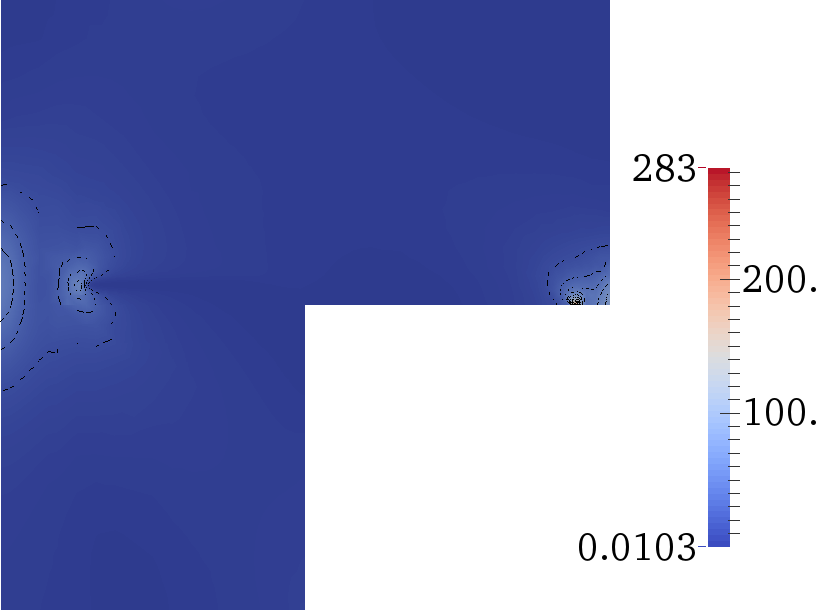}\label{fig:lshape_stresse}} \hspace{10mm}
    \subfloat[]{\includegraphics[scale = 0.35]{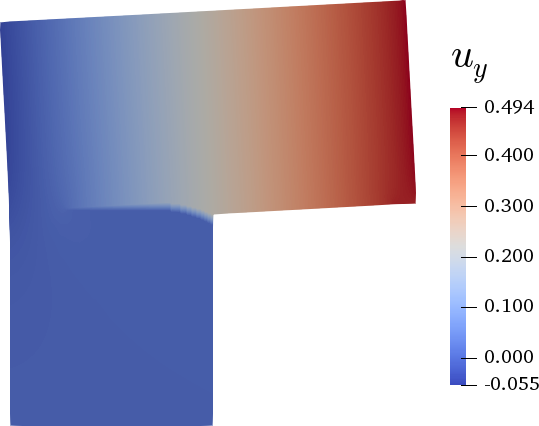}\label{fig:lshape_stressf}}
    \caption{L-shape panel: von-Mises stress distribution and deformed configuration at remote displacement (a, b) 0.23 mm (c, d) 0.29 mm (e, f) 0.49 mm [The deformation is scaled by a factor 50 for visualization].}
    \label{fig:stress_lshape}
\end{figure}

\begin{figure}[H]
    \centering
    \includegraphics[scale = 0.8]{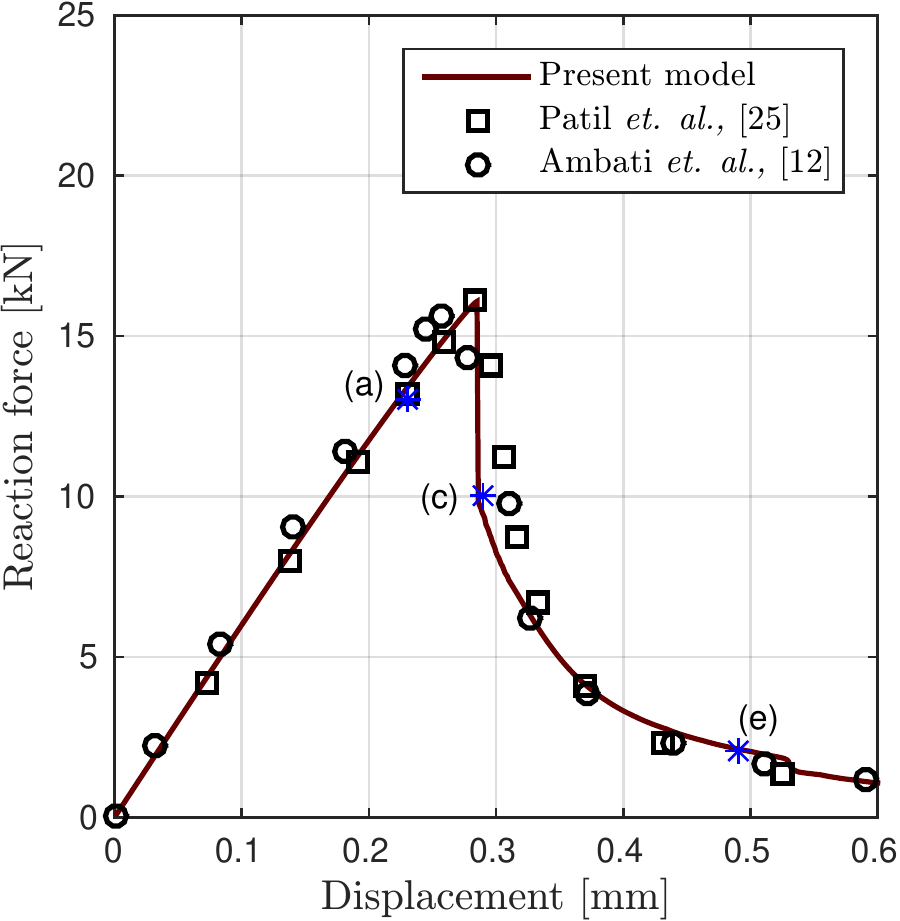}
    \caption{L-shape panel: load-displacement plot [Points A, C and E are the locations of phase field parameter, stress and displacement plots].}
    \label{fig:load_disp_lshape}
\end{figure}

\fref{fig:stress_lshape} shows the von-Mises stress distribution and the displacement map at the different locations (A, C, E) marked on the \fref{fig:load_disp_lshape}. The von-Mises stress concentrates around the corner and the loading region, see \fref{fig:lshape_stressa}. And as the stress at the corner reaches critical, the crack initiates. The location of the von-Mises stress concentration changes as the crack propagates as shown in \frefs{fig:lshape_stressc} and \ref{fig:lshape_stresse}. The deformed configuration shows the discontinuity at the crack location, see \frefs{fig:lshape_stressd} and \ref{fig:lshape_stressf}. \fref{fig:load_disp_lshape} shows the load-displacement response comparison with the standard PFM. The obtained results show a very good agreement with the standard PFM literature \cite{Ambati2015} and adaptive PFM \cite{PATIL2018254}.  

\section{Concluding remarks}
\label{Section:conclusion}
In this paper, an adaptive phase field model is proposed to alleviate the mesh burden associated with the conventional approach. The adaptive refinement technique is based on the recovery based error indicator and the quadtree decomposition. The salient features of the proposed framework are: (a) an efficient recovery based error indicator to identify the elements to be refined; (b) a robust adaptive mesh refinement using quadtrees and (c) elements with hanging nodes are treated as $n-$sided polygonal elements with mean value coordinates as basis functions. With the proposed framework, the potential regions where the fracture can propagate is computed posteriori at each step and the polygons in this region are refined. The performance of the method is systematically demonstrated by solving a few problems in the PFM literature. The results were obtained with significantly fewer degrees of freedom when compared with the the standard PFM results yet demonstrate an excellent agreement. 

\section*{References}
\bibliographystyle{model1-num-names}
\bibliography{myxfemRef,MyRef}

\end{document}